\documentclass[10pt,a4paper,fleqn]{article}

\usepackage[T1]{fontenc}
\usepackage[margin=0.75in]{geometry}
\usepackage{amsmath,amssymb,amsthm}
\usepackage{graphicx}
\usepackage{booktabs}
\usepackage[numbers,sort&compress]{natbib}
\usepackage{placeins}
\usepackage{subcaption}
\usepackage{xcolor}
\usepackage[hidelinks]{hyperref}
\graphicspath{{./}}

\newcommand{\pa}{{\partial}}

\newcommand{\ap}{\alpha}
\newcommand{\be}{\beta}

\newcommand{\util}{\Tilde{u}}
\newcommand{\wtil}{\Tilde{w}}
\newcommand{\vtil}{\Tilde{v}}

\newcommand{\bbu}{\boldsymbol{u}}
\newcommand{\bbg}{\boldsymbol{g}}
\newcommand{\bbtau}{\boldsymbol{\tau}}

\newcommand{\bbn}{\boldsymbol{n}}
\newcommand{\bbU}{\boldsymbol{U}}

\newcommand{\iRe}{\mathrm{Re}^{-1}}

\newcommand{\md}{{\,\rm d}}

\newcommand{\fl}[2]{\frac{#1}{#2}}

\newcommand{\meter}{\mathrm{m}}
\newcommand{\second}{\mathrm{s}}
\newcommand{\kg}{\mathrm{kg}}
\newcommand{\per}{\mathbin{\!/\!}}

\newenvironment{revisionblock}{}{}

\newtheorem{theorem}{Theorem}
\newtheorem{remark}{Remark}
\newtheorem{example}{Example}[section]

\newtheorem{corollary}[theorem]{Corollary}
\newtheorem{proposition}{Proposition}
\newtheorem{definition}{Definition}

\begin{document}

\title{Moment-enhanced shallow-water equations with an effective wall closure for no-slip bottoms}
\author{
\href{https://orcid.org/0000-0002-3247-0328}{Shiping Zhou}$^{1,*}$
\href{https://orcid.org/0000-0003-0527-7431}{Juntao Huang}$^{2}$
\href{https://orcid.org/0000-0002-5395-5455}{Andrew J. Christlieb}$^{1}$\\[0.6em]
\small $^{1}$Computational Mathematics, Science and Engineering, Michigan State University, East Lansing, MI 48824, USA\\
\small $^{2}$Department of Mathematical Sciences, University of Delaware, Newark, DE 19716, USA\\[0.4em]
\small $^{*}$Corresponding author: \href{mailto:zhouship@msu.edu}{zhouship@msu.edu}
}
\date{}

\maketitle

\begin{abstract}
Shallow-water equations and low-order shallow-water moment models use vertically coarse representations and therefore cannot, in general, resolve the thin wall-affected region produced by a no-slip bottom.
Enforcing the pointwise wall value on a low-order global polynomial reconstruction can introduce stiff relaxation and distort the resolved interior velocity profile.
Starting from the incompressible Navier--Stokes equations with Navier bottom friction, we derive a bottom-to-mean relation in a distinguished regular-friction regime and use it to define an endpoint-consistent effective wall-traction closure for the shallow-water equations and the hyperbolic shallow-water moment equations.
The closure represents the momentum effect of unresolved near-wall dynamics; it neither resolves the physical boundary layer nor imposes the pointwise no-slip trace on the reconstructed polynomial.
It recovers the perfect-slip wall contribution when the friction coefficient vanishes.
Because only source terms are changed, the homogeneous principal matrices and their established two-dimensional hyperbolicity classification remain unchanged.
We compare the standard and modified reduced models with two-phase incompressible Navier--Stokes computations in OpenFOAM for wet-bed dam-break and three-dimensional collapse tests.
In the cases considered, the modified closure reduces the excessive damping of the classical low-order wall source and improves agreement in depth-averaged and resolved-interior velocity diagnostics, but it does not uniformly improve front-propagation speed.
The regular-friction asymptotic remainder is not uniform in the large-friction numerical regime; there the effective coefficient is used as a wall-model continuation and assessed empirically.
\end{abstract}


\medskip

\section{Introduction}
Shallow water models have been widely used to describe free-surface flows, including water flow in ocean tides, wave breaking on shallow beaches, roll waves in open channels, flood flows in rivers, surges, dam-break waves, and tsunami propagation.
In addition to hydrodynamic applications, these models are also utilized in studying the dispersion of heavy gases over the Earth's surface and atmospheric dynamics.
If dissipative effects are neglected, then the resulting shallow water equations (SWE) are nonlinear hyperbolic PDEs.
The SWE are based on three main assumptions: (1) a hydrostatic pressure distribution, assuming that vertical acceleration has a negligible effect on pressure; (2) a uniform horizontal velocity distribution along the vertical axis; and (3) water depth is small relative to the wavelength or the curvature of the free surface \cite{Toro2024}.
Compared to the physical complexity and computational cost of the original models, such as the three-dimensional Navier--Stokes equations, shallow water models offer a more efficient and accessible alternative for simulation and analysis.
However, due to bottom friction effects, whereby water near the bottom moves more slowly than water near the surface, the assumption of a constant { vertical profile of the horizontal velocity} leads to inaccuracies.
As a result, shallow water models fail to fully capture variations in the { vertical profile of the horizontal velocity}.

To recover components of the { vertical profile of the horizontal velocity}, one widely studied approach in the literature involves multilayer methods, which decompose the fluid into two or more layers along the vertical direction \cite{Audusse2005, Audusse2008, Fernandez-Nieto2023}.
In \cite{Audusse2005}, the layers are assumed to be advected by the flow without mass exchange between neighboring layers, making the model physically closer to a simulation of non-miscible fluids.
Later, in \cite{Audusse2011}, the model was extended to allow mass and momentum exchanges between layers.
In this formulation, each layer is pre-defined by its thickness and a vertically constant horizontal velocity.
An alternative strategy is the shallow water moment equations (SWME), proposed in \cite{Kowalski2019}, where the { vertical profile of the horizontal velocity} is represented by a polynomial expansion around the mean value.
The expansion coefficients, referred to as moments, allow for the approximation of complex { vertical profiles of the horizontal velocity}.
{ Increasing the retained order enriches the polynomial approximation space, but it does not guarantee uniform improvement of every diagnostic and, for the unregularized SWME, it can destroy hyperbolicity.
In one horizontal dimension, Koellermeier and Rominger introduced the regularized HSWME and $\beta$-HSWME systems \cite{Koellermeier2020}.
In two horizontal dimensions, both first-moment components $(\alpha_1,\beta_1)$ and all higher moment pairs $(\alpha_i,\beta_i)$ enter the state.
The direct two-dimensional HSWME extension retains rotational invariance, but its directional matrix is only weakly hyperbolic at generic states with nonzero first moments.
Bauerle et al. analyzed this distinction and proposed a separate modification that is globally strongly hyperbolic \cite{Bauerle2025}.
The present paper uses the standard direct two-dimensional HSWME principal matrices, not that later globally strongly hyperbolic modification; accordingly, the model is weakly hyperbolic at generic nonzero first-moment states and strongly hyperbolic only when $(\alpha_1,\beta_1)=(0,0)$, as stated in Theorem~\ref{thm:HSWME}.}
Another approach for ensuring hyperbolicity is the shallow water linearized moment equations (SWLME), introduced in \cite{Koellermeier2022}, which neglect nonlinear terms in the higher-order moment equations.
Their stability was analyzed in \cite{Huang2022d}, and well-balanced schemes were developed in \cite{Pimentel-Garcia2024, Caballero-Cardenas2025}.
More recently, several follow-up studies have explored applications and enhancements, including axisymmetric models \cite{Verbiest2025}, and applications to sediment transport \cite{Garres-Diaz2021}.

{ At a solid boundary, the incompressible Navier--Stokes equations are commonly closed by the no-slip condition.
In small-viscosity or high-Reynolds-number regimes, this pointwise constraint produces a thin region of rapid wall-normal velocity variation.
Prandtl's boundary-layer theory explains how a layer that is geometrically thin can nevertheless exert a leading influence on wall shear, momentum transfer, drag, separation, and the outer flow \cite{Prandtl1905,Schlichting2017}.
Direct resolution requires wall-normal resolution tied to the layer thickness and is generally incompatible with the computational purpose of depth-averaged, low-order moment, or otherwise vertically coarse reduced models.}

{ A classical response is to replace unresolved near-wall dynamics by an effective boundary relation or wall law that transmits the required momentum flux to the resolved variables.
This principle underlies wall-function and wall-layer methods in coarse simulations \cite{Launder1974,Piomelli1989,Cabot2000,Piomelli2002,Bose2018} and asymptotically or homogenization-derived effective boundary conditions for flows over unresolved wall structure \cite{Achdou1998,Jager2001}.
In shallow-water asymptotics, related derivations track how viscosity and bottom traction enter depth-averaged equations \cite{Gerbeau2001,Marche2007}.
The purpose of such a closure is not to reconstruct every pointwise feature of the boundary layer on a coarse representation; it is to represent its leading effect on the resolved flow.}

{ This distinction is essential for the present moment formulation.
The horizontal velocity is represented by a low-order global polynomial in the normalized vertical coordinate.
Such a representation cannot, in general, resolve a thin no-slip layer and simultaneously approximate the interior velocity accurately.
In the classical moment source, the large-friction limit acts directly on the reconstructed endpoint value.
This is consistent with the pointwise wall condition, but it imposes one linear constraint on the global polynomial coefficients.
At low order, that localized endpoint constraint can distort the resolved interior profile; increasing the retained order can reduce the profile error, but it does not guarantee improved agreement in all bulk diagnostics.}

{ We therefore treat finite Navier wall friction and its no-slip limit as an effective wall-closure problem for the retained moments.
The underlying Navier--Stokes benchmark imposes the physical wall condition and resolves a wall-affected region to the available mesh scale.
The reduced model instead represents the corresponding wall traction through a bottom-to-mean relation.
It does not claim to resolve the physical boundary layer or to force the reconstructed interior polynomial through the exact wall trace.
Consequently, water height, depth-averaged velocity, resolved-interior vertical profile, and front-propagation speed are distinct diagnostics and are reported separately.}

{ Internal viscous diffusion and bottom friction also have distinct roles.
In this paper, $\mu$ controls viscous momentum diffusion within the fluid, whereas $\kappa$ controls tangential traction in the Navier wall law.
Increasing $\mu$ at fixed slip length tends to relax vertical variation but does not impose zero tangential wall velocity.
No slip is approached through the vanishing dimensionless slip length $\ell_s/H=\iRe/\gamma\to0$.
Under the distinguished regular-friction scaling used in the derivation, the proposed source has an asymptotic wall-reconstruction remainder.
The large-$\kappa$ water calculations lie outside that uniform regime; there the rational coefficient is used as a wall-model continuation and is assessed empirically against the Navier--Stokes computations.}

{ The first contribution of this work is an endpoint-consistent effective wall-traction closure for the SWE and for the standard HSWME principal part specified in Section~\ref{sec:M-SWE}.
At zeroth order, this closure defines the modified shallow water equations (MSWE); combined with the HSWME principal part, it defines the modified hyperbolic shallow water moment equations (MHSWME).
The closure is derived from the bottom-to-mean relation in the regular-friction asymptotic regime and is then used, outside that regime, as a clearly identified wall-model continuation.
Because the modification changes only source terms, the MHSWME retains the homogeneous principal matrices and the two-dimensional hyperbolicity classification of the standard HSWME.
Extension to another friction law, another moment regularization, or a non-flat-bottom regime beyond the retained mild-slope approximation requires a separate derivation and validation.}

{ The second contribution is a systematic numerical comparison of the standard and modified reduced models with two-phase incompressible Navier--Stokes computations using the volume-of-fluid method \cite{Hirt1981} in OpenFOAM \cite{Greenshields2024}.
Direct Navier--Stokes comparisons are not new in viscous shallow-water modeling \cite{Gerbeau2001}; rather, the present study focuses on how the classical and modified low-order moment wall sources affect bulk evolution, resolved-interior vertical profiles, and propagation speed under a large-friction approximation to no slip.
For the test problems considered here, the modified source reduces the excessive damping produced by the classical low-order wall source and improves agreement in the resolved-interior velocity profiles, but it is not uniformly superior in every diagnostic.}

The paper is organized as follows.
{ In Section~\ref{sec:SWE}, we derive the bottom-to-mean relation and its effective wall-traction closure from the dimensional incompressible Navier--Stokes equations with Navier bottom friction and obtain the vertically rescaled formulation.
In Section~\ref{sec:M-SWE}, we introduce the moment expansion and derive the modified hyperbolic shallow water moment equations.}
In Section \ref{sec:numerical}, we conduct numerical simulations and compare their performance.
In Section \ref{sec:conclusion}, we conclude the paper and discuss future directions.
{ \ref{app:other-friction-regimes} discusses other wall-friction regimes, and \ref{app:proof-theorem1} provides the proof of the HSWME hyperbolicity result.}
In addition, to improve the reproducibility of the numerical results, the implementations, including the OpenFOAM setup for solving the incompressible Navier--Stokes equations and our Python code for solving several types of shallow water (moment) equations, are provided in the accompanying source code \cite{Zhou2026}.

\section{{ Asymptotic derivation of the bottom-to-mean relation}}
\label{sec:SWE}

{ In this section, we extend the flat-bottom asymptotic expansion of \cite{Gerbeau2001} to two horizontal dimensions and then to a smooth stationary non-flat bottom.
Starting from the dimensional incompressible Navier--Stokes equations with Navier bottom friction, we apply a shallow-flow nondimensionalization, with $\varepsilon$ denoting the shallowness parameter, and introduce the distinguished scaling of viscosity and wall friction used in the asymptotic expansion.
Under this scaling, the order-$\varepsilon$ correction to the leading depth-uniform velocity has a parabolic vertical profile and yields a bottom-to-mean relation that remains unchanged for a non-flat bottom under the mild-slope assumption through the retained order.
Finally, we derive the vertically rescaled system and the effective wall-traction closure.
For smooth solutions in the regular-friction regime, the corresponding approximation of the exact bottom-traction source has an $\mathcal{O}(\varepsilon^2)$ remainder; other friction scalings are discussed in \ref{app:other-friction-regimes}.}

{ Throughout, $t\geq0$ denotes time, $(x,y)$ are horizontal coordinates, and $z$ is the vertical coordinate.
For a scalar $f$ and a horizontal vector $\boldsymbol{a}=(a_1,a_2)^T$, we use
\begin{equation*}
\nabla_H f:=\left(\fl{\pa f}{\pa x},\fl{\pa f}{\pa y}\right)^T,
\qquad
\nabla_H\!\cdot\boldsymbol{a}:=\fl{\pa a_1}{\pa x}+\fl{\pa a_2}{\pa y},
\qquad
\Delta_H f:=\fl{\pa^2f}{\pa x^2}+\fl{\pa^2f}{\pa y^2},
\end{equation*}
and $\nabla_H^2f$ denotes the horizontal Hessian.
The symbol $I_d$ denotes the $d\times d$ identity matrix, $\boldsymbol{a}\otimes\boldsymbol{b}$ denotes the dyadic product, $|\cdot|$ denotes absolute value for scalars and the Euclidean norm for vectors, and $\|\cdot\|_F$ denotes the Frobenius norm for matrices.
Unless stated otherwise, derivatives are taken at fixed physical coordinates.}

{ Let $\bbu=(\boldsymbol{u}_H,w)^T$, with $\boldsymbol{u}_H=(u,v)^T$, and let the fluid domain be
\begin{equation}\label{eq:dimensional-fluid-domain-corrected}
    \Omega(t)=\left\{(x,y,z):h_b(x,y)<z<h_s(t,x,y)\right\},
    \qquad h=h_s-h_b.
\end{equation}
The dimensional incompressible Navier--Stokes equations are}
\begin{equation}\label{eq:NS-3d}
\left\{
\begin{aligned}
    \nabla\cdot \bbu &= 0, \\
    \fl{\pa\bbu}{\pa t} + \nabla\cdot(\bbu\bbu) &= -\fl{1}{\rho}\nabla p + \fl{1}{\rho}\nabla\cdot\bbtau + \bbg,
\end{aligned}
\right.
\end{equation}
{ where $\bbtau=\mu\big(\nabla\bbu+(\nabla\bbu)^T\big)$ is the viscous stress tensor.}
The density $\rho$ is a constant and $\bbg = (0,0,-g)^{T}$ is the vector of gravitational acceleration with $g>0$.

{ Define the total Cauchy stress tensor as}
\begin{equation}
    { \boldsymbol{T} := -pI_3+\bbtau=-pI_3+\mu\big(\nabla\bbu+(\nabla\bbu)^T\big)},
\end{equation}
where $\mu$ denotes the dynamic viscosity.
The free surface is material and stress-free.
Its kinematic and dynamic boundary conditions are
\begin{equation}\label{bc:freesurface}
    \fl{\pa h_s}{\pa t}+\boldsymbol{u}_H\!\cdot\nabla_Hh_s=w, \quad
    { \boldsymbol{T}\boldsymbol{n}_{s}=\boldsymbol{0}},
\qquad \mbox{ on } z=h_{s}(t,x,y),
\end{equation}
where $h_{s} = h_{s}(t,x,y)$ is the profile of the free surface and $\bbn_{s}$ denotes the outward normal to the free surface.
{ The pressure is measured relative to atmospheric pressure; hence the stress-free condition specifies zero gauge normal pressure together with zero tangential traction at the free surface.}

{
At the smooth stationary bottom, non-penetration and the Navier slip condition \cite{Gerbeau2001} are
\begin{equation}\label{bc:bottom-navier-corrected}
    \bbu\!\cdot\boldsymbol{n}_b=0,
    \quad
    P_b\bbtau\boldsymbol{n}_b=-\kappa P_b\bbu,
    \qquad \mbox{ on } z=h_b,
\end{equation}
where $\boldsymbol{n}_s$ and $\boldsymbol{n}_b$ are the outward unit normals at the free surface and bottom, respectively, $P_b=I_3-\boldsymbol{n}_b\otimes\boldsymbol{n}_b$ is the bottom tangential projector, and $\kappa\geq0$ is the dimensional wall-friction coefficient.
}

\subsection{{ Dimensionless Navier--Stokes system and boundary condition scaling}}
\label{sec:dimensionless-system}

{
Let $L$, $H$, and $U$ be the characteristic horizontal length, vertical length, and horizontal velocity, respectively, and define the shallowness parameter $\varepsilon=H/L$.
We use the shallow-flow scaling
\begin{equation}\label{eq:shallow-scaling-corrected}
\begin{gathered}
    (x,y)=L(\widehat{x},\widehat{y}),
    \qquad z=H\widehat{z}=\varepsilon L\widehat{z},
    \qquad h=H\widehat{h},
    \qquad t=\fl{L}{U}\widehat{t},\\
    \boldsymbol{u}_H=U\widehat{\boldsymbol{u}}_H,
    \qquad w=\varepsilon U\widehat{w},
    \qquad p=\rho U^2\widehat{p}.
\end{gathered}
\end{equation}
The dimensionless groups are
\begin{equation}\label{eq:dimensionless-groups-corrected}
    \iRe=\fl{\mu}{\rho U H},
    \qquad
    G=\fl{gH}{U^2},
    \qquad
    \gamma=\fl{\kappa}{\rho U}.
\end{equation}
{ The parameters $\iRe$ and $\gamma$ describe distinct physical effects: $\iRe$ multiplies viscous diffusion in the bulk equations, while $\gamma$ measures the tangential wall traction in the Navier boundary condition.
The dimensional Navier slip length is $\ell_s=\mu/\kappa$, and hence its dimensionless counterpart is $\ell_s/H=\iRe/\gamma$.
Thus, perfect slip corresponds to $\gamma=0$, whereas the no-slip limit of the Navier slip boundary condition is $\ell_s/H=\iRe/\gamma\to0$ (equivalently, $\gamma/\iRe\to\infty$), for which the tangential wall velocity tends to zero.
Consequently, increasing internal viscosity at a fixed slip condition is not equivalent to imposing no slip: the former relaxes the vertical profile through bulk diffusion, while the latter constrains its value at the wall.}
After dropping the hats, the dimensionless bulk equations are
\begin{equation}\label{sys:NS-dimless-corrected}
\left\{
\begin{aligned}
    \nabla_H\!\cdot\boldsymbol{u}_H+\fl{\pa w}{\pa z}&=0,\\
    \fl{\pa\boldsymbol{u}_H}{\pa t}
    +(\boldsymbol{u}_H\!\cdot\nabla_H)\boldsymbol{u}_H
    +w\fl{\pa\boldsymbol{u}_H}{\pa z}
    +\nabla_Hp
    &=\varepsilon\iRe\Delta_H\boldsymbol{u}_H
    +\fl{\iRe}{\varepsilon}\fl{\pa^2\boldsymbol{u}_H}{\pa z^2},\\
    \varepsilon^2\left[
    \fl{\pa w}{\pa t}
    +\boldsymbol{u}_H\!\cdot\nabla_Hw
    +w\fl{\pa w}{\pa z}
    \right]
    +\fl{\pa p}{\pa z}
    &=-G+\varepsilon\iRe\fl{\pa^2w}{\pa z^2}
    +\varepsilon^3\iRe\Delta_Hw.
\end{aligned}
\right.
\end{equation}
To state the dimensionless boundary conditions without suppressing the bottom geometry, define the scaled velocity and viscous stress
\begin{equation}\label{eq:scaled-velocity-stress-corrected}
    \bbu^{\varepsilon}
    :=\begin{pmatrix}\boldsymbol{u}_H\\ \varepsilon w\end{pmatrix},
    \qquad
    \bbtau^{\varepsilon}
    :=\begin{pmatrix}
    \varepsilon\iRe\left(\nabla_H\boldsymbol{u}_H+(\nabla_H\boldsymbol{u}_H)^T\right)
    &\iRe\left(\fl{\pa\boldsymbol{u}_H}{\pa z}+\varepsilon^2\nabla_Hw\right)\\
    \iRe\left(\fl{\pa\boldsymbol{u}_H}{\pa z}+\varepsilon^2\nabla_Hw\right)^T
    &2\varepsilon\iRe\fl{\pa w}{\pa z}
    \end{pmatrix}.
\end{equation}
The outward unit normals and the bottom tangential projector are
\begin{equation*}
    \boldsymbol{n}_s^{\varepsilon}
    =\fl{(-\varepsilon\nabla_Hh_s,1)^T}
    {\sqrt{1+\varepsilon^2|\nabla_Hh_s|^2}},
    \qquad
    \boldsymbol{n}_b^{\varepsilon}
    =\fl{(\varepsilon\nabla_Hh_b,-1)^T}
    {\sqrt{1+\varepsilon^2|\nabla_Hh_b|^2}},
    \qquad
    P_b^{\varepsilon}=I_3-\boldsymbol{n}_b^{\varepsilon}\otimes\boldsymbol{n}_b^{\varepsilon}.
\end{equation*}
Accordingly, the exact dimensionless free-surface and stationary-bottom conditions are
\begin{align*}
    &\fl{\pa h_s}{\pa t}+\boldsymbol{u}_H\!\cdot\nabla_Hh_s=w, \quad
    \left(-pI_3+\bbtau^{\varepsilon}\right)\boldsymbol{n}_s^{\varepsilon}=\boldsymbol{0}, \qquad \mbox{ on } z=h_s, \\
    &w=\boldsymbol{u}_H\!\cdot\nabla_Hh_b, \quad
    P_b^{\varepsilon}\bbtau^{\varepsilon}\boldsymbol{n}_b^{\varepsilon} =-\gamma P_b^{\varepsilon}\bbu^{\varepsilon}, \qquad \mbox{ on } z=h_b.
\end{align*}
The pressure does not appear in the bottom tangential law because $P_b^{\varepsilon}(-pI_3)\boldsymbol{n}_b^{\varepsilon}=\boldsymbol{0}$.
Together with compatible initial data and appropriate lateral boundary conditions, these equations define the dimensionless free-boundary problem.

For a flat stationary bottom, $\nabla_Hh_b=\boldsymbol{0}$, and the bottom conditions reduce to
\begin{equation}\label{eq:flat-bottom-navier-corrected}
    w=0, \quad
    \gamma\boldsymbol{u}_H -\iRe\left( \fl{\pa\boldsymbol{u}_H}{\pa z} +\varepsilon^2\nabla_Hw \right)=\boldsymbol{0},
    \qquad \mbox{ on } z=h_b.
\end{equation}
The free-surface tangential traction gives
\begin{equation}\label{bc:freesurface-shear-corrected}
    \iRe\fl{\pa\boldsymbol{u}_H}{\pa z}=\mathcal{O}(\varepsilon^2),
    \qquad z=h_s,
\end{equation}
for the expansion below.
Equations \eqref{sys:NS-dimless-corrected}--\eqref{bc:freesurface-shear-corrected} determine the order matching directly in terms of $\iRe$.
}

\subsection{{ Flat-bottom expansion and the bottom-to-mean relation}}
\label{sec:flat-asymptotic-correction}

{
For a flat bottom, set $\eta=z-h_b$, so $0<\eta<h$.
Under
\begin{equation}\label{eq:regular-friction-regime}
    \iRe=\mathcal{O}(1),
    \qquad
    \gamma=\varepsilon\gamma_0,
    \qquad
    \gamma_0=\mathcal{O}(1),
\end{equation}
take the regular expansions
\begin{equation}\label{eq:regular-expansions-manuscript}
\begin{aligned}
    \boldsymbol{u}_H^\varepsilon
    &=\boldsymbol{U}_0+\varepsilon\boldsymbol{u}_1
      +\varepsilon^2\boldsymbol{u}_2+\mathcal{O}(\varepsilon^3),\\
    w^\varepsilon
    &=w_0+\varepsilon w_1+\varepsilon^2w_2+\mathcal{O}(\varepsilon^3),\\
    p^\varepsilon
    &=p_0+\varepsilon p_1+\varepsilon^2p_2+\mathcal{O}(\varepsilon^3).
\end{aligned}
\end{equation}
The $\mathcal{O}(\varepsilon^{-1})$ horizontal equation and the leading bottom and free-surface shear conditions give
\begin{equation}\label{eq:leading-motion-by-slices-manuscript}
    \iRe\fl{\pa^2\boldsymbol{U}_0}{\pa\eta^2}=\boldsymbol{0},
    \qquad
    \fl{\pa\boldsymbol{U}_0}{\pa\eta}(0)
    =\fl{\pa\boldsymbol{U}_0}{\pa\eta}(h)=\boldsymbol{0},
\end{equation}
hence $\boldsymbol{U}_0=\boldsymbol{U}_0(t,x,y)$ is depth-uniform.
Continuity and the leading kinematic conditions then yield
\begin{equation}\label{eq:leading-flat-mass-manuscript}
    \fl{\pa h}{\pa t}+\nabla_H\!\cdot(h\boldsymbol{U}_0)=0.
\end{equation}
The leading pressure is hydrostatic, $p_0=G(h-\eta)$.
For convenience, define the leading horizontal momentum residual
\begin{equation}\label{eq:F0-manuscript}
    \boldsymbol{F}_0
    :=\fl{\pa\boldsymbol{U}_0}{\pa t}
    +(\boldsymbol{U}_0\!\cdot\nabla_H)\boldsymbol{U}_0
    +G\nabla_Hh.
\end{equation}
Substituting the expansions \eqref{eq:regular-expansions-manuscript} into the horizontal momentum equation in \eqref{sys:NS-dimless-corrected}, the flat-bottom Navier condition \eqref{eq:flat-bottom-navier-corrected}, and the free-surface shear condition \eqref{bc:freesurface-shear-corrected}, and then collecting the $\mathcal{O}(1)$ terms, yields the following vertical boundary-value problem for the first-order correction $\boldsymbol{u}_1$ to the horizontal velocity:
\begin{equation}\label{eq:u1-bvp-manuscript}
    \iRe\fl{\pa^2\boldsymbol{u}_1}{\pa\eta^2}=\boldsymbol{F}_0,
    \qquad
    \iRe\fl{\pa\boldsymbol{u}_1}{\pa\eta}(0)=\gamma_0\boldsymbol{U}_0,
    \qquad
    \iRe\fl{\pa\boldsymbol{u}_1}{\pa\eta}(h)=\boldsymbol{0}.
\end{equation}
Depth integration of \eqref{eq:u1-bvp-manuscript} gives $h\boldsymbol{F}_0=-\gamma_0\boldsymbol{U}_0$ and therefore
\begin{equation}\label{eq:leading-flat-momentum-manuscript}
    \fl{\pa(h\boldsymbol{U}_0)}{\pa t}
    +\nabla_H\!\cdot(h\boldsymbol{U}_0\otimes\boldsymbol{U}_0)
    +\fl{G}{2}\nabla_H(h^2)
    =-\gamma_0\boldsymbol{U}_0.
\end{equation}
Solving the local vertical problem and writing $\boldsymbol{u}_b^\varepsilon=\boldsymbol{u}_H^\varepsilon(0)$ gives
\begin{equation}\label{eq:asymptotic-profile-physical}
    \boldsymbol{u}_H^\varepsilon(\eta)
    =\left[
    1+\fl{\gamma}{\iRe}
    \left(\eta-\fl{\eta^2}{2h}\right)
    \right]\boldsymbol{u}_b^\varepsilon
    +\mathcal{O}(\varepsilon^2).
\end{equation}
Consequently, the exact depth average $\boldsymbol{u}_m^\varepsilon=h^{-1}\int_0^h\boldsymbol{u}_H^\varepsilon\md\eta$ satisfies
\begin{equation}\label{eq:wall-mean-relation-physical}
    \boldsymbol{u}_m^\varepsilon
    =\left(1+\fl{h\gamma}{3\iRe}\right)
    \boldsymbol{u}_b^\varepsilon
    +\mathcal{O}(\varepsilon^2).
\end{equation}
This bottom-to-mean relation will be used in Section~\ref{sec:vertical-normalized} to close the bottom-traction term after vertical integration.
}

\subsection{{ Smooth stationary non-flat bottom}}
\label{sec:nonflat-asymptotic-correction}

{
Let $h_b=h_b(x,y)$ be stationary and smooth, and assume the mild-slope geometry
\begin{equation}\label{eq:nonflat-mild-slope-manuscript}
    \varepsilon|\nabla_Hh_b|\ll1,
    \qquad
    |h_b|+|\nabla_Hh_b|+\|\nabla_H^2h_b\|_F=\mathcal{O}(1).
\end{equation}
{ These estimates are understood pointwise and uniformly over the horizontal domain as $\varepsilon\to0$.}
Unlike the flat-bottom specialization \eqref{eq:flat-bottom-navier-corrected}, the exact bottom conditions retain the wall geometry:
\begin{equation}\label{eq:nonflat-bottom-conditions-manuscript}
    w=\boldsymbol{u}_H\!\cdot\nabla_Hh_b,
    \qquad
    P_b^{\varepsilon}\bbtau^{\varepsilon}\boldsymbol{n}_b^{\varepsilon}
    =-\gamma P_b^{\varepsilon}\bbu^{\varepsilon},
    \qquad \mbox{ on } z=h_b.
\end{equation}
Set
\begin{equation*}
    \boldsymbol{b}:=\nabla_Hh_b=(b_1,b_2)^T,
    \qquad
    J_b:=\sqrt{1+\varepsilon^2|\boldsymbol{b}|^2},
    \qquad
    \boldsymbol{a}:=\fl{\pa\boldsymbol{u}_H}{\pa z}+\varepsilon^2\nabla_Hw,
    \qquad
    D_H(\boldsymbol{u}_H):=\nabla_H\boldsymbol{u}_H+(\nabla_H\boldsymbol{u}_H)^T.
\end{equation*}
The two vectors
\begin{equation}\label{eq:nonflat-bottom-tangent-vectors-manuscript}
    { \boldsymbol{t}_{b,1}^{\varepsilon}}
    =\left(1,0,\varepsilon b_1\right)^T,
    \qquad
    { \boldsymbol{t}_{b,2}^{\varepsilon}}
    =\left(0,1,\varepsilon b_2\right)^T,
\end{equation}
span the tangent plane.
Since a vector has zero tangential projection if and only if its scalar product with both tangent vectors vanishes, the projected Navier condition is equivalent to testing $\bbtau^{\varepsilon}\boldsymbol{n}_b^{\varepsilon}+\gamma\bbu^{\varepsilon}$ against { $\boldsymbol{t}_{b,1}^{\varepsilon}$ and $\boldsymbol{t}_{b,2}^{\varepsilon}$}.
Substitution of the scaled stress tensor from \eqref{eq:scaled-velocity-stress-corrected} into the projected Navier condition \eqref{eq:nonflat-bottom-conditions-manuscript}, followed by testing against the tangent vectors \eqref{eq:nonflat-bottom-tangent-vectors-manuscript} and multiplication by the common normalizing factor, gives the horizontal-vector identity
\begin{equation*}
\begin{aligned}
    \gamma J_b\left(\boldsymbol{u}_H+\varepsilon^2w\boldsymbol{b}\right)
    -\iRe\boldsymbol{a}
    +\varepsilon^2\iRe\Big[
        D_H(\boldsymbol{u}_H)\boldsymbol{b}
        +\boldsymbol{b}\left(
        \boldsymbol{a}\!\cdot\boldsymbol{b}
        -2\fl{\pa w}{\pa z}
        \right)
    \Big]
    =\boldsymbol{0},
    \qquad \mbox{ on } z=h_b.
\end{aligned}
\end{equation*}

Introduce the bottom-fitted coordinate $\xi=z-h_b(x,y)$ and use the regular expansions \eqref{eq:regular-expansions-manuscript}, with $\eta$ replaced by $\xi$.
In addition to the mild-slope assumption \eqref{eq:nonflat-mild-slope-manuscript}, we assume the regularity underlying the expansion \eqref{eq:regular-expansions-manuscript}: the expansion coefficients and the horizontal and vertical derivatives appearing below remain uniformly $\mathcal{O}(1)$ as $\varepsilon\to0$.
Under \eqref{eq:nonflat-mild-slope-manuscript},
\begin{equation*}
\begin{aligned}
    J_b&=1+\mathcal{O}(\varepsilon^2),\\
    \boldsymbol{a}\big|_{\xi=0}
    &=\fl{\pa\boldsymbol{U}_0}{\pa\xi}(0)
      +\varepsilon\fl{\pa\boldsymbol{u}_1}{\pa\xi}(0)
      +\mathcal{O}(\varepsilon^2),\\
    \varepsilon^2\iRe\Big[&
        D_H(\boldsymbol{u}_H)\boldsymbol{b}
        +\boldsymbol{b}\left(
        \boldsymbol{a}\!\cdot\boldsymbol{b}
        -2\fl{\pa w}{\pa z}
        \right)
    \Big]=\mathcal{O}(\varepsilon^2).
\end{aligned}
\end{equation*}
Moreover, $\gamma J_b\varepsilon^2w\boldsymbol{b}=\mathcal{O}(\varepsilon^3)$ because $\gamma=\varepsilon\gamma_0$.
Therefore, the $\mathcal{O}(1)$ and $\mathcal{O}(\varepsilon)$ parts of the exact bottom law are, respectively,
\begin{equation*}
    \iRe\fl{\pa\boldsymbol{U}_0}{\pa\xi}(0)=\boldsymbol{0},
    \qquad
    \iRe\fl{\pa\boldsymbol{u}_1}{\pa\xi}(0)=\gamma_0\boldsymbol{U}_0.
\end{equation*}
The free-surface shear condition \eqref{bc:freesurface-shear-corrected} similarly gives
\begin{equation*}
    \iRe\fl{\pa\boldsymbol{U}_0}{\pa\xi}(h)=\boldsymbol{0},
    \qquad
    \iRe\fl{\pa\boldsymbol{u}_1}{\pa\xi}(h)=\boldsymbol{0}.
\end{equation*}
These are precisely the four shear conditions in \eqref{eq:leading-motion-by-slices-manuscript} and \eqref{eq:u1-bvp-manuscript}.

The geometry does change the kinematics and the horizontal forcing.
Define the vertical velocity relative to the bottom by
\begin{equation}\label{eq:relative-vertical-velocity-manuscript}
    q=w-\boldsymbol{u}_H\!\cdot\nabla_Hh_b.
\end{equation}
Then incompressibility becomes $\nabla_H\!\cdot\boldsymbol{u}_H+\pa q/\pa\xi=0$, with $q=0$ at $\xi=0$.
The leading horizontal velocity remains depth-uniform, whereas
\begin{equation}\label{eq:nonflat-leading-kinematics-manuscript}
    q_0=-\xi\nabla_H\!\cdot\boldsymbol{U}_0,
    \qquad
    w_0=\boldsymbol{U}_0\!\cdot\nabla_Hh_b
    -\xi\nabla_H\!\cdot\boldsymbol{U}_0.
\end{equation}
The free-surface kinematic condition again gives \eqref{eq:leading-flat-mass-manuscript}.
The hydrostatic pressure is $p_0=G(h-\xi)=G(h_s-z)$.
Because the horizontal derivative in the momentum equation is taken at fixed physical $z$, its leading gradient is
\begin{equation}\label{eq:nonflat-pressure-gradient-manuscript}
    \nabla_Hp_0\big|_z=G\nabla_H(h+h_b),
\end{equation}
rather than $G\nabla_Hh$ as for a flat bottom.
Consequently, if
\begin{equation*}
    \boldsymbol{F}_{0,b}
    :=\fl{\pa\boldsymbol{U}_0}{\pa t}
      +(\boldsymbol{U}_0\!\cdot\nabla_H)\boldsymbol{U}_0
      +G\nabla_H(h+h_b),
\end{equation*}
then the local first-correction problem is
\begin{equation}\label{eq:nonflat-u1-bvp-manuscript}
    \iRe\fl{\pa^2\boldsymbol{u}_1}{\pa\xi^2}=\boldsymbol{F}_{0,b},
    \qquad
    \iRe\fl{\pa\boldsymbol{u}_1}{\pa\xi}(0)=\gamma_0\boldsymbol{U}_0,
    \qquad
    \iRe\fl{\pa\boldsymbol{u}_1}{\pa\xi}(h)=\boldsymbol{0}.
\end{equation}
It has the same vertical operator and shear data as \eqref{eq:u1-bvp-manuscript}; only the horizontal forcing has changed.
Hence it gives the same parabolic dependence on the distance from the bottom and the same bottom-to-mean relation \eqref{eq:wall-mean-relation-physical} through the retained order.

The assumption $\iRe=\mathcal{O}(1)$ in \eqref{eq:regular-friction-regime} specifies the distinguished asymptotic regime in which the first-order correction and its $\mathcal{O}(\varepsilon^2)$ remainder are derived.
It is not a requirement that the numerical examples must satisfy after the physical parameters have been fixed.
In particular, the no-slip dam-break computations in Section~\ref{sec:numerical} have $\iRe\ll1$ and $\gamma/\varepsilon\gg1$, so the regular expansion is not a uniform error estimate for those computations.
There, the rational effective coefficient obtained from the regular regime is used as a wall-model continuation and is assessed by comparison with the resolved OpenFOAM solutions.
Characteristic reference depth and velocity scales are needed for the asymptotic order interpretation so that the dimensionless variables are order one, whereas a different computational normalization neither changes the dimensional equations nor alters the applicable asymptotic regime.
A consistent change of reference scales does not alter the underlying dimensional problem, but it also cannot, by itself, move that problem into a different asymptotic regime or extend the validity of the remainder estimate.
}

\subsection{Vertically rescaled system}
\label{sec:vertical-normalized}

The shallow water models approximate the incompressible Navier--Stokes equations by averaging over the depth, which results in the loss of information about the { vertical profile of the horizontal velocity}.
To overcome this limitation, the shallow water moment equations (SWME) were proposed by introducing a polynomial expansion in the vertical direction \cite{Kowalski2019}.
In this part, we normalize the dimensionless system in the vertical direction to fully leverage the polynomial basis functions.

{ All mapped-variable derivations below concern wet states with $h(t,x,y)>0$.
Following \cite{Kowalski2019}, define the normalized vertical coordinate
\begin{equation}\label{eq:mapped-coordinate}
    \zeta=\zeta(t,x,y,z):=\fl{z-h_b(x,y)}{h(t,x,y)}\in[0,1].
\end{equation}
Thus $z=h(t,x,y)\zeta+h_b(x,y)$ maps the physical interval $h_b(x,y)\leq z\leq h_s(t,x,y)$ to $0\leq\zeta\leq1$.
For any physical-space field $\theta(t,x,y,z)$, define its mapped counterpart by
\begin{equation}\label{eq:mapped-field-definition}
    \widetilde{\theta}(t,x,y,\zeta)
    :=\theta\bigl(t,x,y,h(t,x,y)\zeta+h_b(x,y)\bigr),
\end{equation}
so that $\theta(t,x,y,z)=\widetilde{\theta}\bigl(t,x,y,\zeta(t,x,y,z)\bigr)$.
The exact mapped bottom traces and depth averages are
\begin{equation}\label{eq:exact-traces-and-means}
\begin{aligned}
u_b^\varepsilon(t,x,y)&:=\widetilde{u}^\varepsilon(t,x,y,0),
&u_m(t,x,y)&:=\int_0^1\widetilde{u}(t,x,y,\widehat\zeta)\md\widehat\zeta,\\
v_b^\varepsilon(t,x,y)&:=\widetilde{v}^\varepsilon(t,x,y,0),
&v_m(t,x,y)&:=\int_0^1\widetilde{v}(t,x,y,\widehat\zeta)\md\widehat\zeta.
\end{aligned}
\end{equation}
The traces $u_b^\varepsilon$ and $v_b^\varepsilon$ are the scalar components of $\boldsymbol{u}_b^\varepsilon$ in \eqref{eq:asymptotic-profile-physical}; they need not equal the endpoint values of a finite Legendre reconstruction.
After \eqref{eq:mapped-field-definition}, explicit dependence on $(t,x,y)$ in $h$ and $h_b$ is suppressed when no ambiguity can arise.
Derivatives of mapped fields are taken at fixed $\zeta$ unless another held-fixed variable is stated.
We also suppress the superscript $\varepsilon$ on bulk mapped fields and their depth averages, while retaining it on the exact wall traces.}

Next we consider the mapping of the mass equation.
Following  \cite{Kowalski2019}, we obtain
\begin{equation}\label{eq:mass-mapped}
    \fl{\pa h}{\pa t} + \fl{\pa (hu_m)}{\pa x} + \fl{\pa (hv_m)}{\pa y} = 0,
\end{equation}
{ where $u_m$ and $v_m$ are the exact mapped depth averages defined in \eqref{eq:exact-traces-and-means}.}

Next, we consider the mapping of the momentum equations.
{ To distinguish the retained shallow-water shear stress from the stress tensor in Section~\ref{sec:dimensionless-system}, define}
\begin{equation}
    { \tau_{xz}  := \fl{\iRe}{\varepsilon}\fl{\pa u}{\pa z}},
\end{equation}
{ whose mapped counterpart is}
\begin{equation}
    { \widetilde{\tau}_{xz}  = \fl{\iRe}{\varepsilon h}\fl{\pa \util}{\pa \zeta}}.
\end{equation}
Then, we have the mapped momentum equation
\begin{equation}\label{eq:u-map1}
\begin{aligned}
    \fl{\pa (h\util)}{\pa t} &+ \fl{\pa}{\pa x}\left( h\util^2 + \fl{G}{2}h^2 \right) 
    + \fl{\pa (h\util\vtil)}{\pa y} \\
    &+ \fl{\pa}{\pa\zeta}\left(  \util\left[\wtil- \fl{\pa (h\zeta+h_b)}{\pa t} -\util\fl{\pa (h\zeta+h_b)}{\pa x} - \vtil\fl{\pa (h\zeta+h_b)}{\pa y} \right] - { \widetilde{\tau}_{xz} } \right)
    = -Gh\fl{\pa h_b}{\pa x}.
\end{aligned}
\end{equation}

Finally, integrating \eqref{eq:u-map1} from 0 to 1 and using the Leibniz integral rule, together with { the retained part of} the bottom Navier boundary condition, gives the mapped $x$-momentum balance
\begin{equation*}
    \fl{\pa (h u_m)}{\pa t} + \fl{\pa }{\pa x}\left(h\int_{0}^{1} \util^2\md\zeta + \fl{G}{2}h^2\right) + \fl{\pa}{\pa y}\left(h\int_{0}^{1}\util\vtil\md\zeta\right) = -\fl{\gamma}{\varepsilon}{ u_b^\varepsilon} - Gh\fl{\pa h_b}{\pa x}.
\end{equation*}
The $y$-momentum equation is derived in the same manner.
Specifically, by defining
\begin{equation*}
    { \tau_{yz} :=\fl{\iRe}{\varepsilon}\fl{\pa v}{\pa z}, \qquad
    \widetilde{\tau}_{yz} =\fl{\iRe}{\varepsilon h}\fl{\pa \vtil}{\pa \zeta}},
\end{equation*}
{ Replacing $\util$, $\widetilde{\tau}_{xz} $, and $\pa h_b/\pa x$ in \eqref{eq:u-map1} by $\vtil$, $\widetilde{\tau}_{yz} $, and $\pa h_b/\pa y$, respectively, gives the corresponding mapped $y$-momentum equation.
Integrating over $\zeta\in[0,1]$ yields the vertically normalized depth-integrated balance laws below.}

\begin{equation}\label{sys:normalized}
\left\{
\begin{aligned}
    \fl{\pa h}{\pa t} + \fl{\pa (hu_m)}{\pa x} + \fl{\pa (hv_m)}{\pa y} &= 0, \\
    \fl{\pa (h u_m)}{\pa t} + \fl{\pa }{\pa x}\left(h\int_{0}^{1} \util^2\md\zeta + \fl{G}{2}h^2\right) + \fl{\pa}{\pa y}\left(h\int_{0}^{1}\util\vtil\md\zeta\right) &= -\fl{\gamma}{\varepsilon}{ u_b^\varepsilon} - Gh\fl{\pa h_b}{\pa x}, \\
    \fl{\pa (h v_m)}{\pa t} + \fl{\pa }{\pa x}\left( h\int_{0}^{1}\util\vtil\md\zeta\right) + \fl{\pa}{\pa y}\left(h\int_{0}^{1}\vtil^2\md\zeta + \fl{G}{2}h^2\right) &= -\fl{\gamma}{\varepsilon}{ v_b^\varepsilon} -Gh\fl{\pa h_b}{\pa y}.
\end{aligned}
\right.
\end{equation}

{ The bottom-to-mean relation \eqref{eq:wall-mean-relation-physical}, applied componentwise to the exact mapped wall traces defined above, gives
\begin{equation}\label{eq:modified-gamma}
\begin{aligned}
-\fl{\gamma}{\varepsilon}u_b^\varepsilon
=-\bar\gamma u_m+\mathcal{O}(\varepsilon^2),\qquad
-\fl{\gamma}{\varepsilon}v_b^\varepsilon
=-\bar\gamma v_m+\mathcal{O}(\varepsilon^2),\qquad
\bar\gamma
:=\fl{\gamma}
{\varepsilon\left(1+\fl{h\gamma}{3\iRe}\right)}.
\end{aligned}
\end{equation}
Here, $u_m$ and $v_m$ are the exact mapped depth averages defined above.
The remainder concerns the approximation of the exact bottom-traction source.
}

\section{{ Moment expansion and modified hyperbolic model}}
\label{sec:M-SWE}

\subsection{Shallow water equations with moment expansion}
\label{sec:sub:SWME}

In this subsection, we derive the moment expansion models of the vertically normalized system introduced in Section \ref{sec:vertical-normalized}.
For notation simplicity, we omit the tilde notation for mapped variables in the remainder of the paper, unless otherwise specified.
We then rewrite the moment expansion equations in a matrix-vector form and provide their hyperbolicity in Theorem \ref{thm:HSWME}.

{ Let $N\in\mathbb{N}_0$ be the retained moment order and define the shifted Legendre basis
\begin{equation}\label{eq:shifted-legendre-basis}
    \phi_j(\zeta):=P_j(1-2\zeta),
    \qquad j\in\mathbb{N}_0,
\end{equation}
where $P_j$ is the degree-$j$ Legendre polynomial.
Then
\begin{equation}\label{eq:shifted-legendre-properties}
    \phi_0\equiv1,
    \qquad \phi_j(0)=1,
    \qquad
    \int_0^1\phi_i(\zeta)\phi_j(\zeta)\md\zeta
    =\fl{\delta_{ij}}{2j+1}.
\end{equation}
The order-$N$ reconstructed mapped velocities are
\begin{align}
\label{eq:u-vertical}
u_N(t,x,y,\zeta)
&=u_m(t,x,y)+\sum_{j=1}^N\alpha_j(t,x,y)\phi_j(\zeta),\\
\label{eq:v-vertical}
v_N(t,x,y,\zeta)
&=v_m(t,x,y)+\sum_{j=1}^N\beta_j(t,x,y)\phi_j(\zeta),
\end{align}
with projection coefficients
\begin{equation}\label{eq:moment-projection-coefficients}
\begin{aligned}
\alpha_j&=(2j+1)\int_0^1\bigl(u-u_m\bigr)\phi_j\md\zeta,\\
\beta_j&=(2j+1)\int_0^1\bigl(v-v_m\bigr)\phi_j\md\zeta,
\qquad j=1,\ldots,N.
\end{aligned}
\end{equation}
Increasing $N$ enriches the polynomial approximation space but does not guarantee uniform improvement of every diagnostic.
The reconstructed bottom values are
\begin{equation}\label{eq:reconstructed-bottom-traces}
\begin{aligned}
u_{b,N}(t,x,y)&:=u_N(t,x,y,0)=u_m(t,x,y)+\sum_{j=1}^N\alpha_j(t,x,y),\\
v_{b,N}(t,x,y)&:=v_N(t,x,y,0)=v_m(t,x,y)+\sum_{j=1}^N\beta_j(t,x,y).
\end{aligned}
\end{equation}
For finite $N$, $u_{b,N}$ and $v_{b,N}$ approximate, but need not equal, the exact physical wall traces $u_b^\varepsilon$ and $v_b^\varepsilon$.
For $N=0$, every sum over $j=1,\ldots,N$ is empty, $u_0=u_m$, and $v_0=v_m$.}

Using the properties of the scaled Legendre polynomials, the shallow water moment equations (SWME) can be derived as \cite{Kowalski2019} (dimensionless):

\begin{align}
\label{eq:SWME-1}    
    \fl{\pa h}{\pa t} &+ \fl{\pa (hu_m)}{\pa x} + \fl{\pa (hv_m)}{\pa y} = 0, \\
\notag    \fl{\pa (hu_m)}{\pa t} &+ \fl{\pa }{\pa x}\left( h(u_m^2+\sum_{j=1}^{N}\fl{\ap_j^2}{2j+1}) + \fl{G}{2}h^2\right) + \fl{\pa}{\pa y}\left( h(u_mv_m + \sum_{j=1}^{N}\fl{\ap_j\beta_j}{2j+1}\right) \\
    &= { -\fl{\gamma}{\varepsilon}u_{b,N}}-Gh\fl{\pa h_b}{\pa x}, \\
\notag    \fl{\pa (hv_m)}{\pa t} &+ \fl{\pa}{\pa x}\left( h(u_mv_m +\sum_{j=1}^{N}\fl{\ap_j\beta_j}{2j+1})\right) + \fl{\pa}{\pa y}\left(h(v_m^2+\sum_{j=1}^{N}\fl{\beta_j^2}{2j+1}) + \fl{G}{2}h^2\right) \\
    &= { -\fl{\gamma}{\varepsilon}v_{b,N}} - Gh\fl{\pa h_b}{\pa y},
\end{align}
with the equations for the higher order averages
\begin{align}
\nonumber    \fl{\pa (h\ap_i)}{\pa t} &+ \fl{\pa}{\pa x}\left( h(2u_m \ap_i + \sum_{j,k=1}^{N}A_{ijk}\ap_j\ap_k)\right) + \fl{\pa}{\pa y}\left(h(u_m\beta_i + v_m\ap_i + \sum_{j,k=1}^{N}A_{ijk}\ap_j\beta_k)\right) \\
    &= u_m D_i - \sum_{j,k=1}^{N} B_{ijk} D_{j}\ap_k { - (2i+1)\left[\fl{\gamma}{\varepsilon}u_{b,N}+\fl{\iRe}{\varepsilon h}\sum_{j=1}^{N} C_{ij}\ap_j\right]}, \\
\nonumber    \fl{\pa (h\beta_i)}{\pa t} &+ \fl{\pa }{\pa x}\left(h(u_m\beta_i + v_m\ap_i +\sum_{j,k=1}^{N}A_{ijk}\ap_j\beta_k)\right) + \fl{\pa}{\pa y}\left( h(2v_m\beta_i + \sum_{j,k=1}^{N}A_{ijk}\beta_j\beta_k)\right) \\
\label{eq:SWME-beta}    &= v_m D_i - \sum_{j,k=1}^{N}B_{ijk}D_j\beta_k { - (2i+1)\left[\fl{\gamma}{\varepsilon}v_{b,N}+\fl{\iRe}{\varepsilon h}\sum_{j=1}^{N} C_{ij}\beta_j\right]},
\end{align}
for $i=1,...,N$, and
\begin{align*}
    A_{ijk} &= (2i+1)\int_{0}^{1}\phi_i \phi_j \phi_k\md\zeta,
    && i,j,k=1,\ldots,N, \\
    B_{ijk} &= (2i+1)\int_{0}^{1}\phi_i'\left(\int_{0}^{\zeta}\phi_j\md\widehat{\zeta}\right)\phi_k\md\zeta,
    && i,j,k=1,\ldots,N,\\
    C_{ij} &= \int_{0}^{1}\phi_i'\phi_j'\md\zeta,
    && i,j=1,\ldots,N,
\end{align*}
where $\phi_i'=\md\phi_i/\md\zeta$.
The coefficient tensors $A_{ijk}$, $B_{ijk}$, and $C_{ij}$ are distinguished from the quasilinear matrices defined below by their indices.
\begin{equation*}
    D_{i} := \fl{\pa (h\ap_i)}{\pa x} +\fl{\pa (h\beta_i)}{\pa y}.
\end{equation*}
This defines the auxiliary quantity $D_i$.

{ Order the state as
\begin{equation}\label{eq:SWME-state-ordering}
    \bbU=(h,hu_m,hv_m,h\alpha_1,h\beta_1,\ldots,h\alpha_N,h\beta_N)^T
    \in\mathbb{R}^{2N+3}.
\end{equation}
After moving every derivative-dependent term in \eqref{eq:SWME-1}--\eqref{eq:SWME-beta} to the left-hand side, define $A^N(\bbU)$ and $B^N(\bbU)$ as the corresponding quasilinear coefficient matrices, with
\begin{equation*}
    A^N(\bbU),B^N(\bbU)\in\mathbb{R}^{(2N+3)\times(2N+3)}.
\end{equation*}
For $N\in\mathbb{N}_0$ and $G>0$, the wet-state admissible set is
\begin{equation}\label{eq:admissible-set}
    \mathcal{A}_N:=\left\{\bbU\in\mathbb{R}^{2N+3}:h>0\right\}.
\end{equation}
The derivations and hyperbolicity statements in this paper are restricted to $\mathcal{A}_N$.
The source vector follows the state ordering in \eqref{eq:SWME-state-ordering}.
\begin{equation}\label{eq:SWME-source-components}
\begin{aligned}
S_{\mathrm{SWME}}^N(\bbU,h_b)=\Bigg[&0,
-\fl{\gamma}{\varepsilon}u_{b,N}-Gh\fl{\pa h_b}{\pa x},
\quad -\fl{\gamma}{\varepsilon}v_{b,N}-Gh\fl{\pa h_b}{\pa y},\\
&-3\left(\fl{\gamma}{\varepsilon}u_{b,N}
+\fl{\iRe}{\varepsilon h}\sum_{j=1}^N C_{1j}\alpha_j\right),
\quad -3\left(\fl{\gamma}{\varepsilon}v_{b,N}
+\fl{\iRe}{\varepsilon h}\sum_{j=1}^N C_{1j}\beta_j\right),\\
&\ldots,\\
&-(2N+1)\left(\fl{\gamma}{\varepsilon}u_{b,N}
+\fl{\iRe}{\varepsilon h}\sum_{j=1}^N C_{Nj}\alpha_j\right),
\quad -(2N+1)\left(\fl{\gamma}{\varepsilon}v_{b,N}
+\fl{\iRe}{\varepsilon h}\sum_{j=1}^N C_{Nj}\beta_j\right)
\Bigg]^T.
\end{aligned}
\end{equation}
For $N=0$, the source vector in \eqref{eq:SWME-source-components} consists only of its first three entries.
The SWME therefore has the self-contained quasilinear form
\begin{equation}\label{sys:SWME-vect}
    \fl{\pa \bbU}{\pa t}
    +A^{N}(\bbU)\fl{\pa \bbU}{\pa x}
    +B^{N}(\bbU)\fl{\pa \bbU}{\pa y}
    =S^{N}_{\mathrm{SWME}}(\bbU,h_b).
\end{equation}}
More details of the matrix-vector form can be found in Section 2 of \cite{Bauerle2025}.

{ The source vector \eqref{eq:SWME-source-components} is nonsingular at $\gamma=0$.
In the perfect-slip case, the wall-traction contributions vanish directly, whereas the vertical-diffusion terms remain.
The resulting source vector is given below.
\begin{equation}\label{eq:SWME-source-perfect-slip}
\begin{aligned}
\left.S_{\mathrm{SWME}}^N(\bbU,h_b)\right|_{\gamma=0}=\Bigg[&0,
-Gh\fl{\pa h_b}{\pa x},
\quad -Gh\fl{\pa h_b}{\pa y},\\
&-3\fl{\iRe}{\varepsilon h}\sum_{j=1}^N C_{1j}\alpha_j,
\quad -3\fl{\iRe}{\varepsilon h}\sum_{j=1}^N C_{1j}\beta_j,\\
&\ldots,\\
&-(2N+1)\fl{\iRe}{\varepsilon h}\sum_{j=1}^N C_{Nj}\alpha_j,
\quad -(2N+1)\fl{\iRe}{\varepsilon h}\sum_{j=1}^N C_{Nj}\beta_j
\Bigg]^T.
\end{aligned}
\end{equation}}

{ The one-dimensional reduction of the SWME system \eqref{sys:SWME-vect} with $N\geq2$ is not globally hyperbolic.
Koellermeier and Rominger \cite{Koellermeier2020} addressed this issue by introducing the hyperbolic shallow water moment equations (HSWME), obtained through a regularization of the one-dimensional system matrix, and established the hyperbolicity of the resulting 1D model.
Bauerle et al.\ \cite{Bauerle2025} subsequently extended this regularized formulation to two spatial dimensions and analyzed its rotational invariance and hyperbolicity.
For $N\geq1$, define the first-moment projection
\begin{equation}\label{eq:first-moment-projection}
    \Pi_1\bbU:=(h,hu_m,hv_m,h\alpha_1,h\beta_1,0,\ldots,0)^T.
\end{equation}
The standard two-dimensional HSWME principal matrices used in this paper are
\begin{equation}\label{eq:HSWME-principal-matrix-definition}
    A_H^N(\bbU):=A^N(\Pi_1\bbU),
    \qquad
    B_H^N(\bbU):=B^N(\Pi_1\bbU),
\end{equation}
and hence $A_H^N,B_H^N\in\mathbb{R}^{(2N+3)\times(2N+3)}$.
The separate globally strongly hyperbolic modification of \cite{Bauerle2025} is not used here.
Following their standard two-dimensional formulation, we write}

\begin{equation}\label{sys:HSWME-vect}
    \fl{\pa \bbU}{\pa t} + A^{N}_{H}(\bbU)\fl{\pa \bbU}{\pa x} + B^{N}_{H}(\bbU)\fl{\pa \bbU}{\pa y} = S^{N}_{\mathrm{SWME}}(\bbU,h_b),
\end{equation}
{ The projection \eqref{eq:first-moment-projection} retains the first-order moments $\alpha_1$ and $\beta_1$ and sets all higher-order moments to zero.
For $N\geq2$, the resulting block patterns are}
\begin{equation}\label{mat:AH}
\begin{aligned}
    A^{N}_{H}(\bbU) = \begin{pmatrix}
0 &1 & & & & & & &  & \\
- u_m^2 - \fl{\ap_1^2}{3} + Gh &2u_m &0 &\fl{2\ap_1}{3} & & & & & & \\
- u_m v_m - \fl{\alpha_1 \beta_1}{3} & v_m & u_m & \frac{\beta_1}{3} & \frac{\alpha_1}{3} & & & & & \\
-2u_m \alpha_1 & 2\alpha_1 & 0 & u_m & \frac{3}{5}\alpha_1 & & & & & \\
-(u_m \beta_1 + v_m \alpha_1) & \beta_1 & \alpha_1 & 0 & u_m & \fl{1}{5}\beta_1 &\fl{2}{5}\ap_1 & & & \\
- \frac{2}{3} \alpha_1^2 & 0 & 0 & \frac{1}{3}\alpha_1 &0 &u_m  &0 & \ddots & & \\
- \frac{2}{3} \alpha_1 \beta_1 & 0  &0 & -\frac{1}{3}\beta_1 &\fl{2}{3}\ap_1 & 0 &u_m & \ddots &0 &0 \\
 &  &  &  &  & \ddots & \ddots & \ddots &\fl{N+1}{2N+1}\ap_1 &0\\
& & & & &\ddots &\ddots &\ddots & \frac{1}{2N+1}\beta_1 & \frac{N}{2N+1}\alpha_1 \\
& & & & &  & \frac{N-1}{2N-1}\alpha_1 & 0 &u_m &0 \\
& & & && &-\fl{1}{2N-1}\beta_1 &\fl{N}{2N-1}\ap_1 &0 & u_m
\end{pmatrix}
\end{aligned},
\end{equation}
and
\begin{equation}\label{mat:BH}
\begin{aligned}
    B^{N}_{H}(\bbU) = \begin{pmatrix}
0 &0 &1 & & & & & &  & \\
- u_m v_m - \fl{\ap_1\beta_1}{3}  &v_m &u_m &\fl{\beta_1}{3} &\fl{\ap_1}{3} & & & & & \\
- v_m^2 - \fl{\beta_1^2}{3} + Gh & 0 & 2v_m & 0 & \frac{2\beta_1}{3} & & & & & \\
-(u_m\beta_1+v_m\ap_1) &\beta_1 & \alpha_1 & v_m &0 & \frac{2}{5}\beta_1 &\fl{1}{5}\ap_1 & & &  \\
-2v_m\beta_1 &0 & 2\beta_1 & 0 & v_m &0 &\fl{3}{5}\beta_1 & & & \\
- \frac{2}{3} \alpha_1\beta_1 & 0 & 0 & \frac{2}{3}\beta_1 &0 &v_m  &0 & \ddots & & \\
- \frac{2}{3} \beta_1^2 & 0  &0 &0 &\fl{1}{3}\beta_1 & 0 &v_m & \ddots &0 &0 \\
 &  &  &  &  & \ddots & \ddots & \ddots &\fl{N+1}{2N+1}\beta_1 &\fl{1}{2N+1}\ap_1\\
& & & & &\ddots &\ddots &\ddots &0 & \frac{N+1}{2N+1}\beta_1 \\
& & & & &  & \frac{N}{2N-1}\beta_1 &-\fl{1}{2N-1}\ap_1 &v_m &0 \\
& & & && &0 &\fl{N-1}{2N-1}\beta_1 &0 & v_m
\end{pmatrix}
\end{aligned}.
\end{equation}

{ Since hyperbolicity is determined by the homogeneous principal matrices $A_H^N$ and $B_H^N$, the source term does not affect the classification.
The two-dimensional hyperbolicity classification is stated in the following theorem:}

\begin{theorem}[hyperbolicity of the HSWME]\label{thm:HSWME}
    { Let $N\geq1$, $G>0$, and $\bbU\in\mathcal{A}_N$.
The HSWME \eqref{sys:HSWME-vect} is weakly hyperbolic in 2D at $\bbU$.
Moreover, it is strongly hyperbolic in 2D at $\bbU$ if and only if $(\alpha_1,\beta_1)=(0,0)$.}
\end{theorem}
\begin{proof}
    { See the proof in \ref{app:proof-theorem1}.}
\end{proof}

\begin{revisionblock}
{ To focus on the friction-induced stiffness, we take $h_b\equiv\mbox{const}$ in the analysis and numerical tests considered below.
For smooth non-flat bottom topography, the standard source terms $-Gh\pa_xh_b$ and $-Gh\pa_yh_b$ remain in the depth-averaged momentum equations in \eqref{sys:normalized}.
However, the accuracy of the wall closure introduced in Section~\ref{sec:M-SWE} also depends on the vertical variation of the local momentum residual and therefore requires separate assessment for non-flat topography.
Consider the zeroth-order system (SWE) with}
\end{revisionblock}

\begin{equation}\label{sys:SWE}
\begin{aligned}
    A_{H}^{0}(\bbU) = \begin{pmatrix}
    0 &1 &0 \\
    -u_m^2+Gh &2u_m &0 \\
    -u_mv_m & v_m &u_m
    \end{pmatrix}, \quad B_{H}^{0}(\bbU) = \begin{pmatrix}
        0 &0 &1 \\
        -u_m v_m &v_m &u_m \\
        -v_m^2+Gh &0 &2v_m
    \end{pmatrix}, \quad S_{\mathrm{SWME}}^{0}(\bbU)=-\fl{\gamma}{\varepsilon}\begin{pmatrix}
        0 \\ u_m \\ v_m
    \end{pmatrix}
\end{aligned},
\end{equation}
and the first-order SWME with
\begin{equation}\label{sys:SWME1}
\begin{aligned}
    A_H^{1}(\bbU) & = \begin{pmatrix}
        0 &1 &0 &0 &0 \\
        -u_m^2+Gh-\fl{1}{3}\ap_1^2 & 2u_m &0 &\fl{2}{3}\ap_1 &0 \\
        -u_m v_m-\fl{1}{3}\ap_1\beta_1 & v_m &u_m &\fl{1}{3}\beta_1 &\fl{1}{3}\ap_1 \\
        -2u_m\ap_1 & 2\ap_1 &0 &u_m &0 \\
        -u_m\beta_1-v_m\ap_1 & \beta_1 &\ap_1 &0 &u_m
    \end{pmatrix}, \\
    B_H^{1}(\bbU) & = \begin{pmatrix}
        0 &0 &1 &0 &0 \\
        -u_mv_m-\fl{1}{3}\ap_1\beta_1 &v_m & u_m &\fl{1}{3}\beta_1 &\fl{1}{3}\ap_1 \\
        -v_m^2+Gh-\fl{1}{3}\beta_1^2 &0 &2v_m &0 &\fl{2}{3}\beta_1 \\
        -u_m\beta_1-v_m\ap_1 &\beta_1 &\ap_1 &v_m &0 \\
        -2v_m\beta_1 &0 &2\beta_1 &0 &v_m
    \end{pmatrix}.
\end{aligned}
\end{equation}
and the source term
\begin{equation}\label{eq:SWME:s1}
    { S_{\mathrm{SWME}}^{1}(\bbU) = \begin{pmatrix}
        0 \\
        -\fl{\gamma}{\varepsilon}(u_m+\ap_1) \\
        -\fl{\gamma}{\varepsilon}(v_m+\beta_1) \\
        -3\left[\fl{\gamma}{\varepsilon}(u_m+\ap_1)+4\fl{\iRe}{\varepsilon h}\ap_1\right] \\
        -3\left[\fl{\gamma}{\varepsilon}(v_m+\beta_1)+4\fl{\iRe}{\varepsilon h}\beta_1\right]
    \end{pmatrix}.}
\end{equation}

{ Equations~\eqref{sys:SWE}--\eqref{eq:SWME:s1} give the order-zero matrices $A_H^0,B_H^0$, the order-one matrices $A_H^1,B_H^1$, and their corresponding sources $S_{\mathrm{SWME}}^0,S_{\mathrm{SWME}}^1$ separately from the general $N\geq2$ block displays.
During a source substep, the water depth is constant and the source operator is linear in the velocity moments.
Hence $\gamma/\varepsilon\gg 1$ produces a rapidly relaxing linear system and imposes an additional stability restriction on a fully explicit discretization.
This restriction can be treated efficiently using an implicit or semi-implicit source update, as in \cite{Huang2022d}.

In the large-friction regime used here to approximate no slip, the standard $N$-moment closure replaces the exact physical wall traces $u_b^{\varepsilon}$ and $v_b^{\varepsilon}$ by the reconstructed endpoint values $u_{b,N}$ and $v_{b,N}$.
The resulting source therefore penalizes the reconstructed bottom velocity.
For $N=0$, $u_{b,0}=u_m$, and $v_{b,0}=v_m$, so the source-relaxation substep directly damps the depth-averaged velocities.
For $N\geq 1$, the stiff relaxation instead drives $u_{b,N}=u_m+\sum_{j=1}^N\alpha_j$ and $v_{b,N}=v_m+\sum_{j=1}^N\beta_j$ toward zero, thereby imposing an endpoint constraint on each global polynomial reconstruction.
Although this treatment is consistent with the pointwise wall condition at the reconstructed level, a finite-order global polynomial may accommodate the endpoint constraint by distorting the resolved interior profile.
Increasing $N$ enlarges the approximation space, but it does not remove the distinction between the exact physical wall trace and the reconstructed endpoint value, nor does it guarantee uniform improvement of every solution diagnostic.

This distinction motivates the effective wall-traction source introduced in the next subsection.
Rather than applying the wall law directly to the reconstructed endpoint values, the modified closure uses the bottom-to-mean relation derived in Section~\ref{sec:SWE} to represent the wall-traction contribution in terms of the retained depth-averaged variables.
The resulting MHSWME is defined for every retained order $N\in\mathbb{N}_0$ and does not impose the exact pointwise wall trace on the truncated Legendre reconstruction.
Its purpose is to represent the leading influence of the unresolved wall traction on the retained variables, not to resolve the physical near-wall boundary layer.
In the regular-friction regime, the underlying bottom-to-mean relation has the stated $\mathcal{O}(\varepsilon^2)$ remainder for smooth solutions; in the large-friction calculations, it is used as a wall-model continuation.
The numerical behavior of the standard and modified closures across retained orders and across distinct diagnostics is assessed separately in Section~\ref{sec:numerical}.
}

\subsection{{ Effective wall-traction source and the MHSWME}}
\label{sec:corrected-MHSWME}

{
At the vertically continuous level, the wall-traction terms contain the exact physical traces $u_b^\varepsilon$ and $v_b^\varepsilon$.
The standard $N$-moment model closes these terms with the reconstructed endpoint values $u_{b,N}$ and $v_{b,N}$, whereas the modified model applies the effective relation \eqref{eq:modified-gamma} directly to the exact wall traction.
This replacement changes only the wall-traction part of the moment source; it does not replace the Legendre reconstruction in the advective fluxes, and it leaves the vertical-diffusion terms containing $C_{ij}$ unchanged.
For the mean equations, the exact wall contributions $-(\gamma/\varepsilon)u_b^\varepsilon$ and $-(\gamma/\varepsilon)v_b^\varepsilon$ are replaced by $-\bar\gamma u_m$ and $-\bar\gamma v_m$, respectively.
In the $i$-th pair of moment equations, the corresponding replacement is
\begin{equation}\label{eq:moment-wall-source-replacement}
    -(2i+1)\fl{\gamma}{\varepsilon}(u_b^\varepsilon,v_b^\varepsilon)^T
    \longmapsto
    -(2i+1)\bar\gamma(u_m,v_m)^T.
\end{equation}
This closure does not impose $(u_{b,N},v_{b,N})=(u_b^\varepsilon,v_b^\varepsilon)$ on the truncated Legendre reconstruction.
Accordingly, the modified hyperbolic shallow water moment equations are
\begin{equation}\label{sys:MHSWME-vect}
    \fl{\pa\bbU}{\pa t}
    +A_H^N(\bbU)\fl{\pa\bbU}{\pa x}
    +B_H^N(\bbU)\fl{\pa\bbU}{\pa y}
    =S_{\mathrm{MHSWME}}^N(\bbU,h_b),
\end{equation}
with the final source formula
\begin{equation}\label{eq:final-MHSWME-source}
\begin{aligned}
S_{\mathrm{MHSWME}}^N(\bbU,h_b)=\Bigg[&0, -\bar\gamma u_m-Gh\fl{\pa h_b}{\pa x},
\quad -\bar\gamma v_m-Gh\fl{\pa h_b}{\pa y},\\
&-3\left(\bar\gamma u_m
+\fl{\iRe}{\varepsilon h}\sum_{j=1}^N C_{1j}\alpha_j\right),
\quad
-3\left(\bar\gamma v_m
+\fl{\iRe}{\varepsilon h}\sum_{j=1}^N C_{1j}\beta_j\right),\\
&\ldots,\\
&-(2N+1)\left(\bar\gamma u_m
+\fl{\iRe}{\varepsilon h}\sum_{j=1}^N C_{Nj}\alpha_j\right), -(2N+1)\left(\bar\gamma v_m
+\fl{\iRe}{\varepsilon h}\sum_{j=1}^N C_{Nj}\beta_j\right)
\Bigg]^T,
\end{aligned}
\end{equation}
where $\bar\gamma$ is given by \eqref{eq:modified-gamma}.
For $N=0$, \eqref{sys:MHSWME-vect} reduces to the modified shallow water equations with source
\begin{equation}\label{eq:final-MSWE-source}
    S_{\mathrm{MSWE}}^0
    =\left(0,
    -\bar\gamma u_m-Gh\fl{\pa h_b}{\pa x},
    -\bar\gamma v_m-Gh\fl{\pa h_b}{\pa y}
    \right)^T.
\end{equation}
For fixed $\varepsilon, h>0$, the algebraic continuation satisfies
\begin{equation}\label{eq:modified-gamma-large-friction-limit}
    \bar\gamma\longrightarrow\fl{3\iRe}{\varepsilon h}
    \qquad\mbox{as}\qquad \gamma\longrightarrow\infty.
\end{equation}
This finite saturation concerns the modeled wall coefficient.
It does not make the regular-friction $\mathcal{O}(\varepsilon^2)$ remainder uniform as $\gamma\to\infty$ and does not resolve a physical no-slip boundary layer.
}

\begin{corollary}[{ Hyperbolicity of the MHSWME}]
{ Let $N\geq1$, $G>0$, and $\bbU\in\mathcal{A}_N$.
The MHSWME \eqref{sys:MHSWME-vect} has the same homogeneous principal matrices as the HSWME.
Therefore, its weak and strong two-dimensional hyperbolicity classification is exactly that stated in Theorem~\ref{thm:HSWME}.}
\end{corollary}

{ The numerical tests below may pass through states with nonzero first moments, for which the standard two-dimensional HSWME is only weakly hyperbolic.
The computations therefore assess the behavior of the stated discretization for these test cases; they do not constitute a general well-posedness result for the weakly hyperbolic system.}

\section{Numerical simulations}
\label{sec:numerical}

In this section, we present numerical simulations to compare the performance of different models with the Navier--Stokes equations.
Here, we use OpenFOAM (VOF) to simulate the incompressible Navier--Stokes equations.
{ All reduced-model tests considered here are wet-bed problems with $h>0$; dry states and wetting--drying interfaces are not treated.}
For comparison, we consider the following models:
\begin{enumerate}
    \item the classical shallow water equation (SWE)
    \item the hyperbolic shallow water moment equations (HSWME)
    \item the modified classical shallow water equation (MSWE)
    \item the modified hyperbolic shallow water moment equations (MHSWME)
\end{enumerate}

{ Across the examples below, these four reduced models are compared with OpenFOAM in terms of water height, depth-averaged horizontal velocity, vertical profiles of the horizontal velocity, and water-front propagation speed.
The 2D dam-break example additionally examines OpenFOAM vertical-resolution and interface-threshold sensitivity, contrasts perfect-slip and no-slip bottoms, and studies the effect of the retained moment order $N$, while the two 3D examples consider zero and linear initial velocity profiles.}

{ The dimensional 2D reference solutions are computed using OpenFOAM's laminar \texttt{incompressibleVoF} two-phase solver.
The computational domain is $[0,100]\times[0,2]\,\meter^2$ in the horizontal--vertical plane, with the water--air interface represented by the water volume fraction $\alpha_{\rm w}$.
The physical parameters for water are $\rho_{\rm w}=1000\,\kg\per\meter^3$, $\nu_{\rm w}=10^{-6}\,\meter^2\per\second$, and for air $\rho_{\rm a}=1\,\kg\per\meter^3$, $\nu_{\rm a}=1.48\times10^{-5}\,\meter^2\per\second$, surface tension $\sigma=0.07\,\mathrm{N}\per\meter$, and gravitational acceleration $g=9.81\,\meter\per\second^2$.
The OpenFOAM calculation and the reduced models do not solve the same continuum problem: OpenFOAM resolves a viscous water--air system with an interface-capturing method and surface tension, whereas the reduced systems are hydrostatic depth-integrated single-fluid models with a modeled bottom traction and no resolved air phase or VOF interface.
The comparisons below are therefore model-to-benchmark comparisons of selected diagnostics, rather than discretization-error comparisons for identical equations.

A uniform Cartesian mesh with $N_x=4000$ and $N_z=160, 320, $ or $640$ is employed to assess sensitivity to the vertical resolution while keeping the horizontal VOF resolution fixed.
The initial time step is $10^{-3}\,\second$ and is subsequently adjusted to ensure that both the flow and interface Courant numbers remain below one.
Unless otherwise stated, the results obtained with $N_z=640$ are used as the OpenFOAM reference solutions.
All simulations are performed up to $t=3\,\second$.
For velocity $\boldsymbol{u}$, the 2D inlet maintains the prescribed profile, the outlet uses zero normal gradient, the upper atmosphere uses \texttt{pressureInletOutletVelocity}, and the lower wall uses either \texttt{slip} or \texttt{noSlip} as specified in Example~\ref{ex:2d-slip-noslip}.
The modified pressure $p_{\rm rgh}$ has zero normal gradient at the inlet, outlet, and lower wall and \texttt{prghTotalPressure} at the atmosphere.
The water fraction $\alpha_{\rm w}$ has zero normal gradient at the inlet, outlet, and lower wall and an air-valued \texttt{inletOutlet} condition at the atmosphere; no contact-angle condition is imposed in these cases.
The bottom elevation is constant, $h_b=0$.
The OpenFOAM \texttt{noSlip} condition imposes the wall velocity exactly in the discrete Navier--Stokes problem.
By contrast, $\kappa=10^{10}\,\kg\per(\meter^2\cdot\second)$ in the reduced models is a finite large-friction approximation to the no-slip limit, and the modified coefficient~\eqref{eq:modified-gamma} is an effective wall model; neither should be identified with an exactly imposed pointwise no-slip trace.

The OpenFOAM time derivative is discretized by Euler's method; the momentum convection uses \texttt{linearUpwind} with a linear gradient, and the VOF transport uses \texttt{interfaceCompression} with \texttt{vanLeer} and two volume-fraction corrections per step.
The pressure--velocity coupling uses three correctors, while linear interpolation and uncorrected normal gradients and Laplacians are used on the orthogonal Cartesian meshes.

From the cell-centered OpenFOAM fields, let $\ap_{ij}$ and $u_{ij}$ denote the water volume fraction and horizontal velocity in the $j$th cell of the vertical column centered at $x_i$, and define $\bar{\ap}_{ij}=\min(1,\max(0,\ap_{ij}))$.
From the uniform vertical spacing $\Delta z$, we compute
\begin{equation}\label{eq:openfoam-diagnostics}
    h_{\rm OF}(x_i,t)=\sum_j\bar\alpha_{ij}\Delta z,
    \qquad
    q_{\rm OF}(x_i,t)=\sum_j\bar\alpha_{ij}u_{ij}\Delta z,
    \qquad
    u_{m,{\rm OF}}(x_i,t)=\frac{q_{\rm OF}(x_i,t)}{h_{\rm OF}(x_i,t)}.
\end{equation}

Thus, partially filled interface cells contribute proportionally to both the water height and the depth-integrated horizontal volume flux (discharge).
To extract a vertical profile at a prescribed horizontal location $x_\star$, we use $\ap_{\rm w}$ and $u$ from the closest cell-center column satisfying $x_i\leq x_\star$, after which only cells satisfying $\ap_{\rm w}\geq \ap_{\star}$ are retained.
Thus, for the uniform mesh with $N_x=4000$, a profile requested at $x_\star=55\,\meter$ uses the column centered at $x_i=54.9875\,\meter$.
Unless otherwise stated, we use $\ap_{\star}=1$.
No horizontal interpolation is applied to these 2D profile data or to the corresponding reduced-model reconstruction: both use the closest cell center at or to the left of $x_\star$.
For the 3D slice and profile comparisons, the nearest available horizontal cell center is used independently in each coordinate, with the lower-coordinate center selected in an exact tie, again without spatial interpolation.
The OpenFOAM vertical profile retains only cells with $\alpha_{\rm w}\geq\alpha_\star$ and is augmented by the exact physical no-slip datum $u=0$ at the bottom; the reduced profile is the Legendre reconstruction from the selected depth-averaged cell and is not augmented by, or fitted to, that datum.

For the front-speed diagnostic, a $0.5\,\meter$ box average is used only to localize strong descending-gradient candidates in $50.5\leq x\leq90\,\meter$.
Candidates have at least $75\%$ of the strongest descending-gradient magnitude.
For each candidate, $h_{\rm behind}$ and $h_{\rm ahead}$ are the medians in the intervals $1.5$--$2.5\,\meter$ behind and ahead of the candidate, and a transition is accepted only if its jump exceeds both $0.02\,\meter$ and five times the robust ahead-state noise estimate.
The front position is the descending crossing of
\begin{equation}\label{eq:front-definition}
    h_f=h_{\rm ahead}+\frac{1}{2}\left(h_{\rm behind}-h_{\rm ahead}\right),
\end{equation}
linearly interpolated between the two adjacent height samples that bracket $h_f$ within $4\,\meter$ of the candidate and recentered once.
Linear interpolation is used only to locate the front crossing; field and vertical-profile comparisons use selected cell-center data without interpolation.
The reported effective speed is the least-squares slope of the accepted front positions over the output times stated in the corresponding caption.

Although the OpenFOAM solutions resolve the full two-phase Navier--Stokes system and therefore provide a higher-fidelity numerical benchmark for comparison, they should not be interpreted as an exact pointwise reference solution.
In particular, VOF-based free-surface solvers can exhibit discretization and interface-capturing artifacts near the water--air interface.
For example, free-surface wiggles and oscillatory or overpredicted velocities have been reported for \texttt{interFoam}, with their magnitude depending on the spatial and temporal resolution, Courant number, and interface treatment \cite{Larsen2019, Ferro2022}.
At the same time, oscillatory structures in a wet-bed dam-break problem are not necessarily numerical artifacts.
For sufficiently large downstream-to-upstream depth ratios, nonhydrostatic and dispersive effects are known to produce an undular bore consisting of a leading wave followed by a train of secondary waves \cite{Kim2011,Cantero-Chinchilla2020}.
In the wet-bed dam-break examples considered below, the depth ratio is $h_R/h_L=1/1.5=2/3$, which lies in the regime in which such undular behavior has been reported.
Therefore, the oscillations observed behind the propagating front may contain physical nonhydrostatic wave structures, whereas small-scale variations localized near the VOF interface may also contain numerical contributions.
Accordingly, we use the OpenFOAM solution primarily as a higher-fidelity Navier--Stokes numerical benchmark for the resolved water height, depth-averaged velocity, vertical profile of the horizontal velocity, and propagation speed, rather than regarding every local interfacial feature as an exact reference value.

The OpenFOAM reference solutions for the remaining numerical examples are computed using the same physical parameters and numerical settings unless otherwise stated.
To avoid repetition, only the problem-specific computational domain and mesh resolution are reported for those examples.

For the reduced models, we use a first-order cell-centered path-conservative finite-volume discretization \cite{Castro2017}.
The nonconservative products are evaluated along straight-line paths using midpoint quadrature at the arithmetic interface states, with local Lax--Friedrichs stabilization and no higher-order reconstruction.
The directional wave-speed estimates are $|u_m|+\sqrt{Gh+\alpha_1^2}$ and $|v_m|+\sqrt{Gh+\beta_1^2}$, with $\alpha_1=\beta_1=0$ for the depth-averaged equations; these estimates are also used for weakly hyperbolic HSWME states.

The transport terms are advanced explicitly by forward Euler.
We use $\Delta t=0.7\Delta x/a_{\max}$ in one horizontal dimension and $\Delta t=0.7/(a_{x,\max}/\Delta x+a_{y,\max}/\Delta y)$ in two horizontal dimensions, where the maximum directional speeds are evaluated over all cells.
For a consistent comparison, the linear sources of all classical and modified moment models are advanced after the transport step using the same local backward-Euler procedure.
The transported depth is frozen during this source substep, and the directional momentum blocks are solved independently by direct linear solves.

All tests use strictly positive wet-bed data; no positivity correction is applied, and the calculation is terminated if a nonpositive depth occurs.
One boundary-cell layer is used on each side.
The 1D and radial-collapse calculations use adjacent-cell zero-normal-gradient copies, with corner states in the latter taken from the nearest diagonally adjacent interior cell.
In Example~\ref{ex:3d-nonzero}, the outflow copies at $x=100\,\meter$ and $y=100\,\meter$ are imposed before the fixed inflows at $x=0$ and $y=0$, which are obtained by exact moment projection of the prescribed linear profile.
}

To improve the reproducibility of the numerical results, the implementations, including the OpenFOAM setup for solving the incompressible Navier--Stokes equations and our Python code for solving several types of shallow water (moment) equations, are provided in the accompanying source code \cite{Zhou2026}.

{
\begin{example}[2D dam-break problem with perfect-slip and no-slip bottoms]\label{ex:2d-slip-noslip}
We consider a two-dimensional dam-break problem on $[0,100]\times[0,2]\,\meter^2$, with horizontal coordinate $x$, vertical coordinate $z$, and flat bottom $h_b=0$.
The initial water depth is
\begin{equation}\label{eq:ex-new-dambreak-height}
    h_0(x)=
    \begin{cases}
        1.5\,\meter, & 0\leq x<50\,\meter,\\
        1.0\,\meter, & 50\,\meter\leq x\leq100\,\meter.
    \end{cases}
\end{equation}
Let $\zeta=z/h_0(x)$.
We compare the following two bottom conditions using initial profiles that have the same depth-averaged horizontal velocity.

For the perfect-slip case, OpenFOAM uses
\begin{equation}\label{eq:ex-new-perfect-slip-bc}
    w=0,
    \qquad
    \frac{\partial u}{\partial z}=0,
    \qquad \mbox{on }z=0,
\end{equation}
which is imposed by the \texttt{slip} condition on the lower wall.
The initial velocity is the compatible cubic profile
\begin{equation}\label{eq:ex-new-perfect-slip-initial}
    u(0,x,z)=A\left(3\zeta^2-2\zeta^3\right),
    \qquad
    w(0,x,z)=0,
    \qquad
    A=\frac{1}{3}\,\meter\per\second.
\end{equation}
It satisfies $\partial_z u=0$ at both $z=0$ and $z=h_0(x)$.
In the shallow water and moment models, perfect slip corresponds to $\kappa=0$.

For the no-slip case, OpenFOAM imposes
\begin{equation}\label{eq:ex-new-no-slip-bc}
    u=w=0,
    \qquad \mbox{on }z=0,
\end{equation}
using the \texttt{noSlip} lower-wall condition.
The initial horizontal velocity is the quadratic profile
\begin{equation}\label{eq:ex-new-no-slip-initial}
    u(0,x,z)=S\left(2\zeta-\zeta^2\right),
    \qquad
    w(0,x,z)=0,
    \qquad
    S=0.25\,\meter\per\second,
\end{equation}
which satisfies $u=0$ at the bottom and $\partial_z u=0$ at the water surface.
In the shallow water and moment models, we take $\kappa=10^{10}\,\kg\per(\meter^2\cdot\second)$ to approximate the no-slip limit.
The two initial profiles are normalized so that
\begin{equation}\label{eq:ex-new-equal-mean}
    \frac{1}{h_0}\int_0^{h_0}u_{\rm slip}(0,x,z)\md z
    =\frac{A}{2}
    =\frac{2S}{3}
    =\frac{1}{6}\,\meter\per\second
    =\frac{1}{h_0}\int_0^{h_0}u_{\rm no\text{-}slip}(0,x,z)\md z.
\end{equation}
Consequently, the comparison isolates the effect of the bottom condition and the associated vertical shape without changing the initial discharge.
Both cases use the material parameters and OpenFOAM settings stated above and are simulated until $t=3\,\second$.
For the reduced models, the computational nondimensionalization uses the reference scales $L=100\,\meter$, $H=1.5\,\meter$, and $U_{\rm ref}=100\,\meter\per\second$.
The essential dimensionless and discretization parameters are listed in Table~\ref{tab:ex-new-dimensionless-parameters}.

\begin{center}
\small
\captionof{table}{{ Parameters used by the reduced models in Example~\ref{ex:2d-slip-noslip}.}}
\label{tab:ex-new-dimensionless-parameters}
\begin{tabular}{@{}lll@{}}
\toprule
{ Quantity} & { Definition} & { Value} \\
\midrule
{ Shallowness parameter}
& $\varepsilon=H/L$
& $0.015$ \\
{ Gravity parameter}
& $G=gH/U_{\rm ref}^{2}$
& $1.4715\times10^{-3}$ \\
{ Inverse Reynolds number}
& $\iRe=\nu/(U_{\rm ref}H)$
& $6.6667\times10^{-9}$ \\
{ Vertical-viscosity coefficient}
& $\iRe/\varepsilon$
& $4.4444\times10^{-7}$ \\
{ Wall-friction parameter}
& $\gamma=\kappa/(\rho U_{\rm ref})$
& $0$ (perfect slip), $10^{5}$ (no-slip approximation) \\
{ Reduced-model discretization}
& $N_x$, CFL
& $4000$, $0.7$ \\
\bottomrule
\end{tabular}
\end{center}

The initial physical flow is subcritical.
Using the common initial depth-averaged velocity $u_m=1/6\,\meter\per\second$ and $1\leq h_0\leq1.5\,\meter$ gives $0.043\leq u_m/\sqrt{gh_0}\leq0.053$; even the Froude number based on the maximum initial pointwise velocity remains below $0.11$.
The reference velocity $U_{\rm ref}$ is used only for computational normalization, so $G^{-1/2}$ does not represent the physical Froude number.

For the dimensionless vertical coordinate $\zeta=z/h_0(x)\in[0,1]$, the scaled Legendre basis functions $\phi_j(\zeta)$ are dimensionless.
The two initial velocity profiles, retained here in dimensional form, have the expansions
\begin{equation}\label{eq:ex-new-initial-moments}
\begin{aligned}
    u_{\rm slip}(0,x,z)
    &=\frac{A}{2}-\frac{3A}{5}\phi_1(\zeta)
      +0\phi_2(\zeta)+\frac{A}{10}\phi_3(\zeta),\\
    u_{\rm no\text{-}slip}(0,x,z)
    &=\frac{2S}{3}-\frac{S}{2}\phi_1(\zeta)
      -\frac{S}{6}\phi_2(\zeta)+0\phi_3(\zeta).
\end{aligned}
\end{equation}
Thus, $N=3$ exactly reconstructs the perfect-slip cubic, while its $N=1$ and $N=2$ projections coincide because the second coefficient vanishes.
The no-slip quadratic is already represented exactly for $N=2$.
Because $\zeta$ and the basis functions $\phi_j(\zeta)$ are dimensionless, all expansion coefficients in \eqref{eq:ex-new-initial-moments} have units of velocity.
For the reduced-model computation, these coefficients are normalized by $U_{\rm ref}$ during initialization and converted back to dimensional units for reporting.
\end{example}
}

{
We first assess the OpenFOAM solutions independently of the shallow water models.
Figures~\ref{fig:ex-new-openfoam-h} and \ref{fig:ex-new-openfoam-um} assess vertical-mesh sensitivity for the no-slip case.
The corresponding perfect-slip bulk diagnostics are omitted because they are nearly indistinguishable at the displayed scale.
Figures~\ref{fig:ex-new-openfoam-u} and \ref{fig:ex-new-openfoam-alpha} then assess the vertical profiles and their sensitivity to the VOF threshold.
The prescribed state at $t=0$ is omitted because it is already given by \eqref{eq:ex-new-dambreak-height}, \eqref{eq:ex-new-perfect-slip-initial}, and \eqref{eq:ex-new-no-slip-initial}.

\begin{figure}[htbp]
\centering
\includegraphics[width=0.96\linewidth]{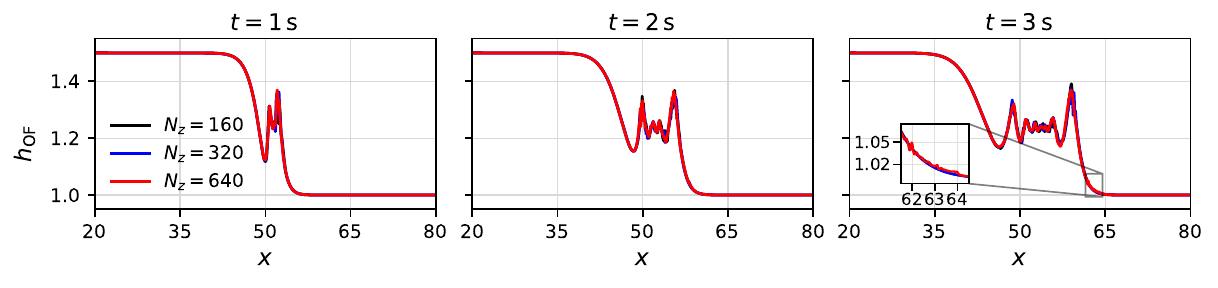}
\caption{{ OpenFOAM vertical-resolution study for the no-slip 2D wet-bed dam-break problem: water height $h_{\rm OF}(t,x)$ at $t=1$, $2$, and $3\,\second$ from left to right, displayed on $20\leq x\leq80\,\meter$.
At fixed $N_x=4000$, the black, blue, and red curves use $N_z=160$, $320$, and $640$, respectively.
The three meshes give nearly coincident bulk profiles and downstream-front locations, with a small spread near the rapidly varying downstream transition at $t=3\,\second$ (inset).
The $N_z=640$ result is used below as the finest available reference, without interpreting the local interfacial structures as pointwise-converged features.}}
\label{fig:ex-new-openfoam-h}
\end{figure}
\FloatBarrier

{ As the disturbed region broadens, its downstream transition propagates to the right and an oscillatory structure develops between the two transitions.
Such behavior is consistent with an undular bore at the depth ratio $h_R/h_L=2/3$, for which nonhydrostatic and dispersive effects can generate a leading wave followed by secondary waves \cite{Kim2011,Cantero-Chinchilla2020}.
The persistence of the large-scale structure across the tested meshes supports its use as a bulk benchmark, whereas the remaining small-scale differences can reflect VOF interface-capturing and discretization effects.
We therefore do not regard the individual interfacial irregularities as exact pointwise reference values.}

\begin{figure}[htbp]
    \centering
    \includegraphics[width=0.96\linewidth]{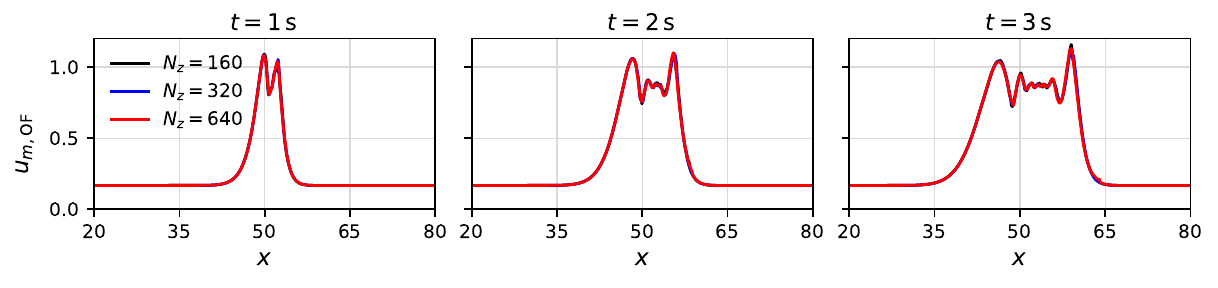}
    \caption{{ OpenFOAM vertical-resolution study for the no-slip 2D wet-bed dam-break problem: depth-averaged horizontal velocity $u_{m,{\rm OF}}(t,x)=q_{\rm OF}/h_{\rm OF}$ at $t=1$, $2$, and $3\,\second$ from left to right, displayed on $20\leq x\leq80\,\meter$.
At fixed $N_x=4000$, the black, blue, and red curves use $N_z=160$, $320$, and $640$, respectively.
The curves are nearly coincident through the bulk flow; the only clearly visible mesh-dependent difference is localized near $x=63$--$64\,\meter$ at $t=3\,\second$ in the rapidly varying downstream transition.
Thus, this bulk velocity diagnostic is weakly sensitive to vertical resolution over the tested range at the displayed scale.}}
    \label{fig:ex-new-openfoam-um}
\end{figure}
\FloatBarrier

{ The localized difference at $t=3\,\second$ occurs where the ratio $u_{m,{\rm OF}}=q_{\rm OF}/h_{\rm OF}$ reflects cell-to-cell variations in both the VOF-weighted discharge and the water height.
This supports using the $N_z=640$ solution for the subsequent bulk comparison, but it does not establish pointwise convergence within the interfacial transition.}

\begin{figure}[htbp]
    \centering
    \includegraphics[width=0.48\linewidth]{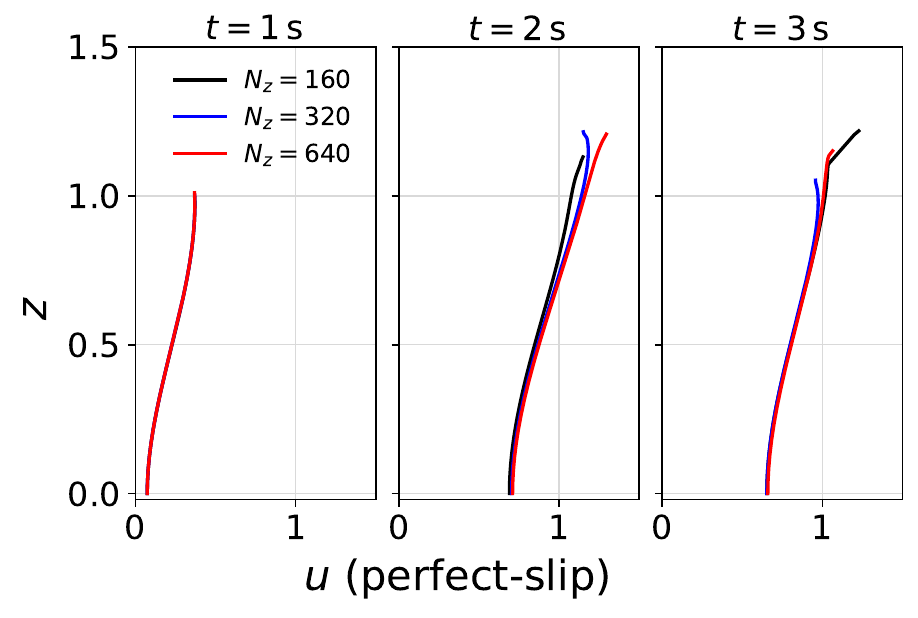}
    \includegraphics[width=0.48\linewidth]{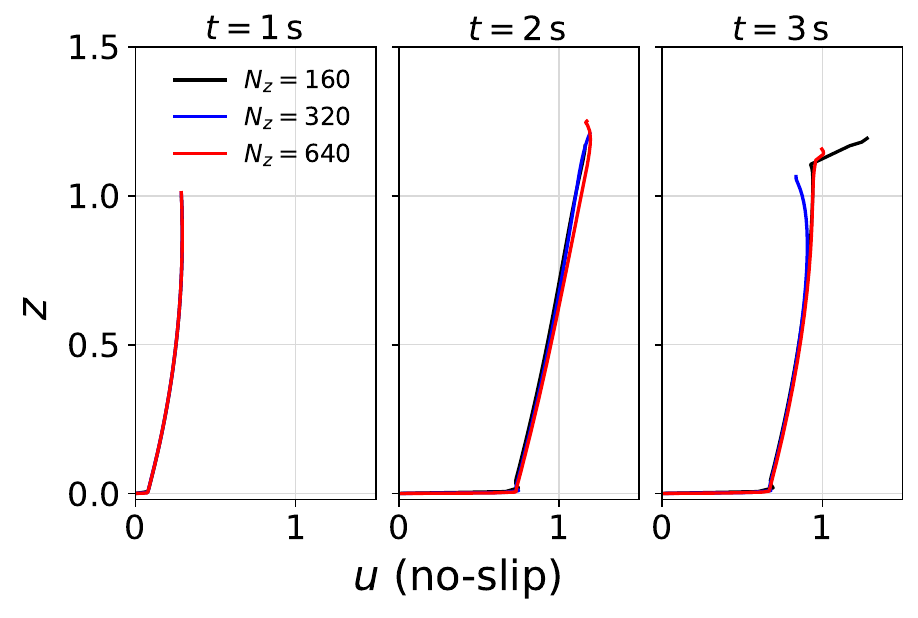}

    \caption{{ OpenFOAM vertical profiles of horizontal velocity $u(z)$ at $x=55\,\meter$ for perfect-slip (left group) and no-slip (right group) bottom conditions.
Within each group, the panels show $t=1$, $2$, and $3\,\second$ from left to right; black, blue, and red denote $N_z=160$, $320$, and $640$ at fixed $N_x=4000$, and only pure-water cell centers satisfying $\alpha_{\rm w}\geq\alpha_\star=1$ are retained.
Perfect slip retains a nonzero tangential wall velocity, whereas no slip connects to $u=0$ through a thin wall-affected region.
The resolutions agree closely through the resolved interior, with their largest differences near the water--air interface.
The perfect-slip wall values are reconstructed from adjacent cell centers, while the exact no-slip datum $u=0$ is added at $z=0$ for visualization.}}
    \label{fig:ex-new-openfoam-u}
\end{figure}
\FloatBarrier

{ The contrast between the two groups isolates the effect of the bottom condition on the vertical structure.
The interior agreement over the tested meshes supports using the $N_z=640$ profiles as the subsequent OpenFOAM benchmark, while the cell-centered wall reconstruction and the threshold-dependent interface region are not used to assess pointwise convergence.}

\FloatBarrier

\begin{figure}[htbp]
    \centering
    \includegraphics[width=0.96\linewidth]{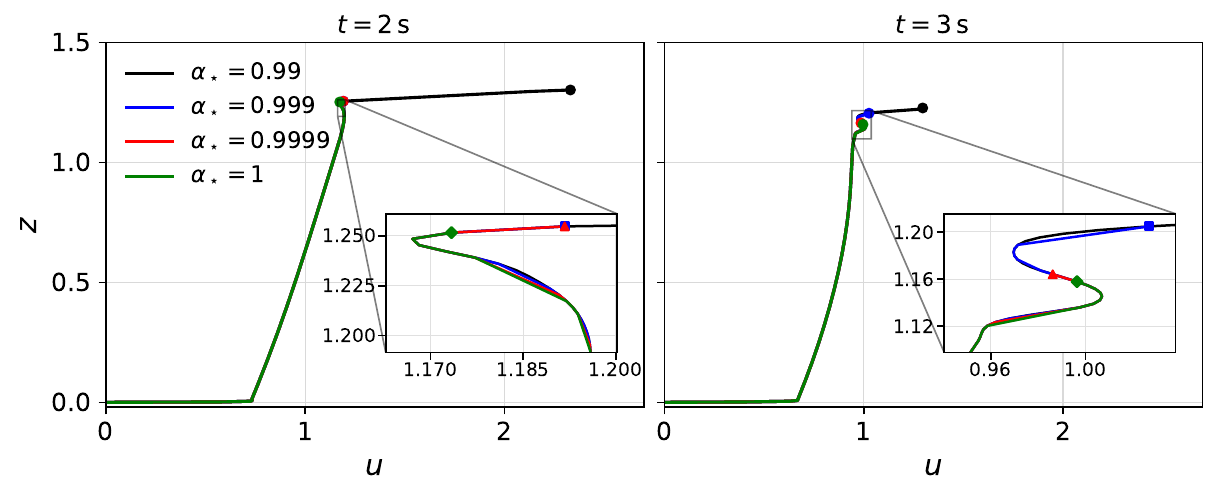}
    \caption{{ Sensitivity of the no-slip OpenFOAM vertical velocity profile at $x=55\,\meter$ to the VOF water-selection threshold $\alpha_\star$ on the $N_x\times N_z=4000\times640$ mesh.
The panels show $t=2$ and $3\,\second$; curves correspond to $\alpha_\star=0.99$, $0.999$, $0.9999$, and $1$, symbols mark the highest retained cell center, and the insets magnify the water--air interface.
The threshold changes only the uppermost mixed cells, while the profiles coincide through the resolved pure-water interior.
Consequently, $\alpha_\star=1$ is used below and threshold-dependent endpoints are not interpreted as pointwise liquid velocities at the free surface.}}
    \label{fig:ex-new-openfoam-alpha}
\end{figure}

\FloatBarrier

{ Lower thresholds retain additional mixed cells and, at $t=2\,\second$, produce a high-speed upper segment that disappears as $\alpha_\star$ approaches one.
Because a mixed cell contains both water and air, its one-fluid velocity need not represent the pointwise liquid velocity.
The retained cells are identical across thresholds at $t=0$ and $1\,\second$, so those times are omitted.}

\FloatBarrier

\begin{figure}[htbp]
    \centering
    \includegraphics[width=0.96\linewidth]{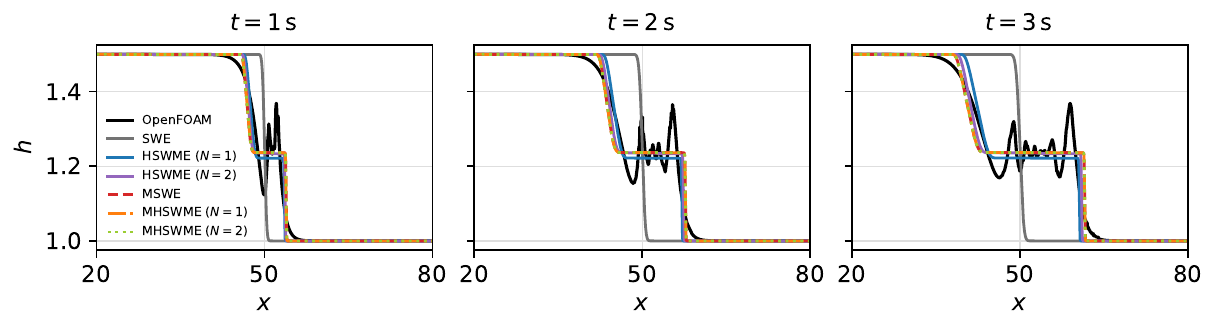}
    \par\vspace{-3mm}
    \includegraphics[width=0.96\linewidth]{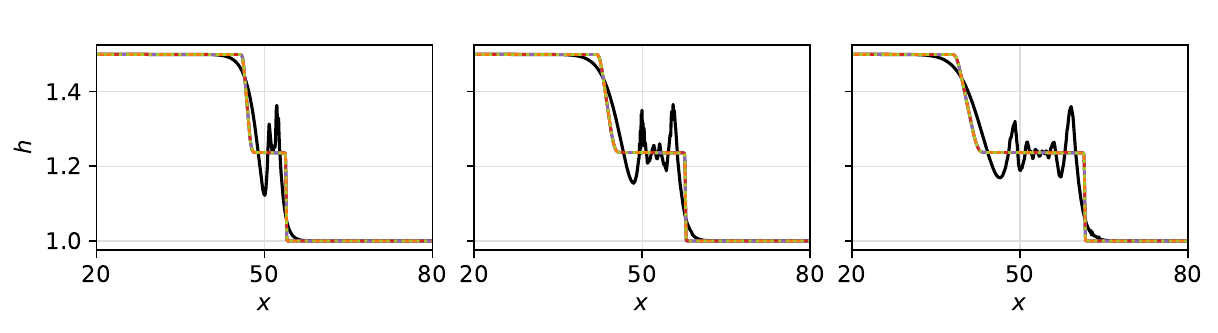}
    \caption{{ Water height $h(t,x)$ in the 2D wet-bed dam-break problem at $t=1$, $2$, and $3\,\second$ from left to right.
The top row compares the no-slip OpenFOAM reference with reduced models using $\kappa=10^{10}\,\kg\per(\meter^2\cdot\second)$; the bottom row uses perfect slip with $\kappa=0$.
Curves identify OpenFOAM, SWE, HSWME with $N=1,2$, MSWE, and MHSWME with $N=1,2$.
Under perfect slip, each original model coincides with its modified counterpart and the reduced models capture the bulk motion of the two transitions, but not the localized OpenFOAM dispersive and VOF structures.
Under no slip, the classical SWE remains close to the initial step and the HSWME fronts are retarded, whereas the modified models recover propagating transitions, although front-location errors remain.}}
    \label{fig:ex-new-model-h}
\end{figure}
\FloatBarrier

{ The overlap of each original/modified pair under perfect slip confirms that the modification is inactive when the wall traction vanishes.
Under no slip, the contrast instead isolates the effect of the wall source: the classical source produces strong damping, while the effective source restores bulk motion.
The water-height curves alone do not establish a uniformly most accurate model, because HSWME with $N=2$ is comparable to the modified models at some displayed times and the front-speed and vertical-profile diagnostics behave differently.}

\begin{figure}[htbp]
    \centering
    \includegraphics[width=0.96\linewidth]{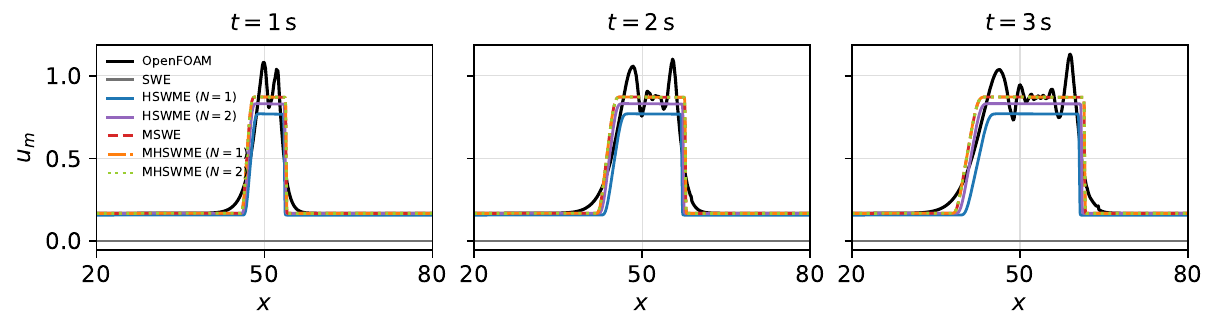}
    \par\vspace{-3mm}
    \includegraphics[width=0.96\linewidth]{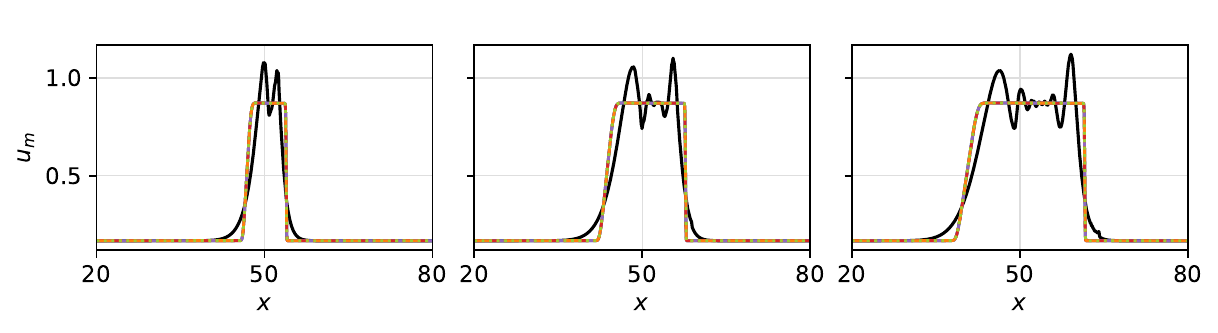}
    \caption{{ Depth-averaged horizontal velocity $u_m(t,x)$ in the 2D wet-bed dam-break problem at $t=1$, $2$, and $3\,\second$ from left to right.
The top row is the no-slip comparison, with $\kappa=10^{10}\,\kg\per(\meter^2\cdot\second)$ in the reduced models, and the bottom row is the perfect-slip comparison with $\kappa=0$.
Curves compare OpenFOAM, SWE, HSWME with $N=1,2$, MSWE, and MHSWME with $N=1,2$.
Under perfect slip, original and modified counterparts overlap and capture the bulk high-velocity region.
Under no slip, the SWE is nearly arrested and the HSWME remains order-dependent, while the modified models follow the bulk OpenFOAM velocity more closely; localized VOF-scale irregularities near the moving front are not expected to be resolved by the depth-averaged systems.}}
    \label{fig:ex-new-model-um}
\end{figure}
\FloatBarrier

{ The mean-velocity comparison reinforces the source-induced damping seen in the water height.
The modified closure recovers the level and extent of the bulk high-velocity region more closely than the classical no-slip source, but it cannot reproduce local OpenFOAM variations that are filtered by depth averaging.
In particular, the irregularity near $x=63$--$64\,\meter$ at $t=3\,\second$ is the localized VOF and horizontal-discretization sensitivity identified in Figure~\ref{fig:ex-new-openfoam-um}, rather than a target feature of the reduced systems.}

\begin{figure}[htbp]
    \centering
    \includegraphics[width=0.80\linewidth]{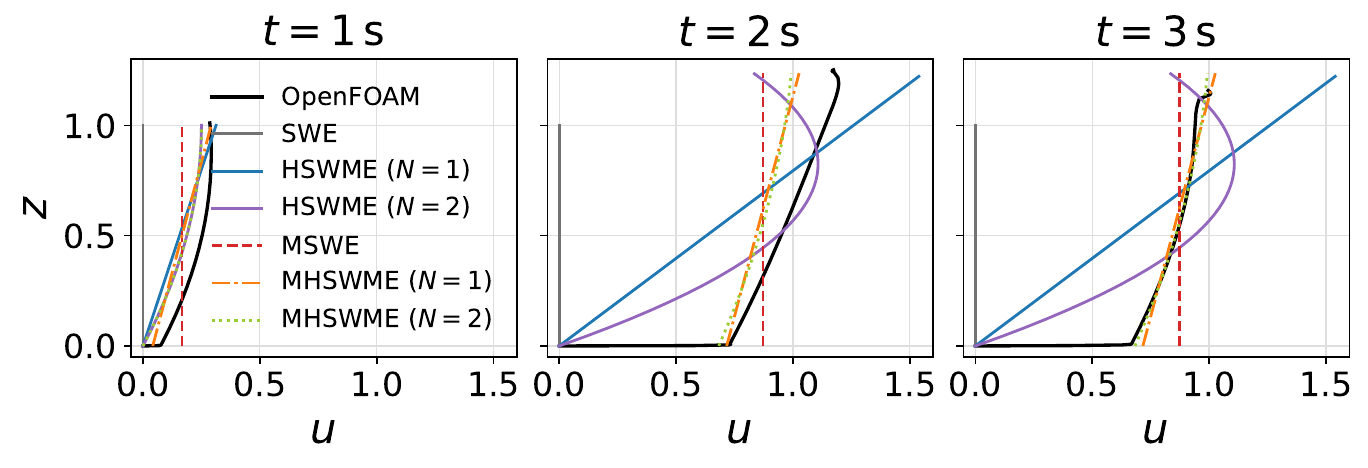}
    \par\vspace{-2mm}
    \includegraphics[width=0.80\linewidth]{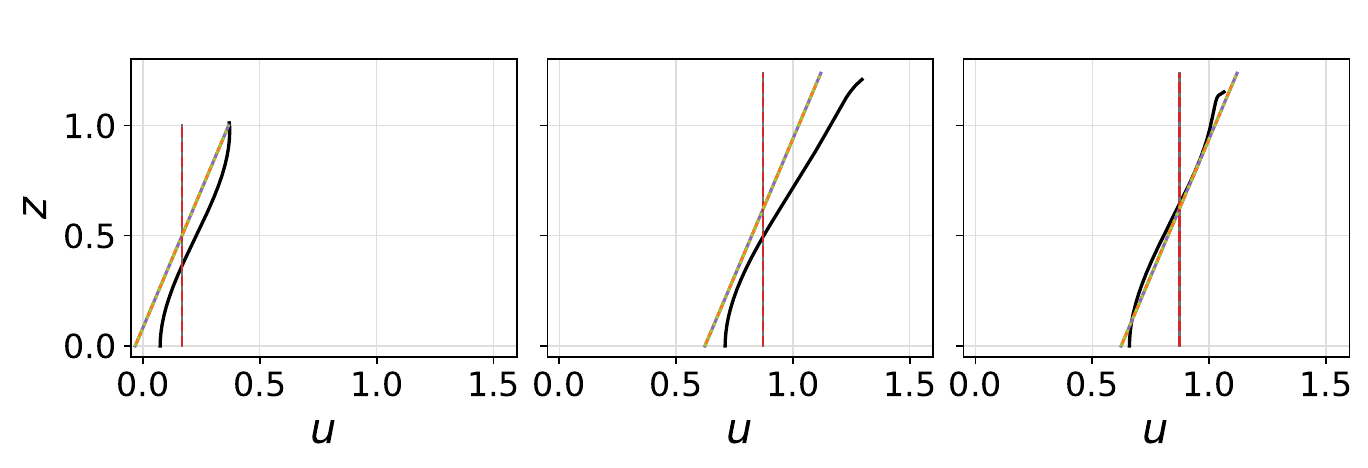}
    \caption{{ Vertical profiles of horizontal velocity $u(z)$ at $x=55\,\meter$ in the 2D wet-bed dam-break problem at $t=1$, $2$, and $3\,\second$.
The top row is no slip and the bottom row is perfect slip.
OpenFOAM profiles retain only pure-water cell centers ($\alpha_\star=1$); HSWME and MHSWME profiles are reconstructed with $N=1$ or $2$.
Under perfect slip, each original model overlaps its modified counterpart and the two retained orders are visually indistinguishable at this location.
Under no slip, the low-order HSWME reconstructions accommodate the zero wall endpoint by strongly distorting the resolved interior profile, whereas the modified models follow the OpenFOAM interior more closely without claiming to resolve the thin physical wall layer.}}
    \label{fig:ex-new-model-u}
\end{figure}
\FloatBarrier

{ Under perfect slip, equation~\eqref{eq:ex-new-initial-moments} gives $\alpha_2(0,x)=0$, and the nonlinear evolution generates only a small second moment at $x=55\,\meter$: its magnitude remains below $1.01\times10^{-5}\,\meter\per\second$ at the displayed times.
Consistently, the domain-wide maxima of $|u_m^{N=2}-u_m^{N=1}|$ are $1.11\times10^{-5}$, $1.92\times10^{-5}$, and $3.00\times10^{-5}\,\meter\per\second$ at $t=1$, $2$, and $3\,\second$, so the two retained orders cannot be separated visually here.
In contrast, the no-slip quadratic has the nonzero initial coefficient $\alpha_2=-S/6=-4.17\times10^{-2}\,\meter\per\second$, and the classical low-order global reconstructions accommodate the stiff endpoint source by departing strongly from the OpenFOAM interior at later times.
The modified source better represents that resolved interior without resolving the physical wall layer.
We therefore assess front speed separately and then test whether increasing $N$ improves the classical HSWME profile.}

\FloatBarrier

\begin{figure}[htbp]
    \centering
    \includegraphics[width=0.96\linewidth]{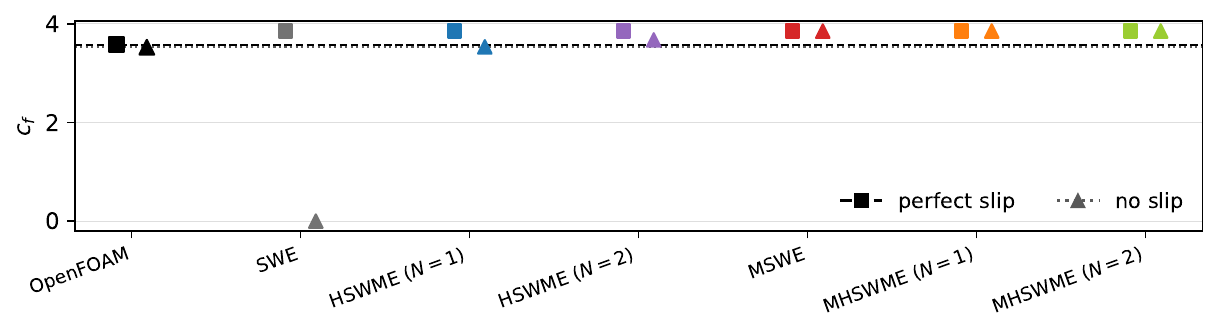}
    \caption{{ Effective right-moving front speed $c_f$ for the 2D wet-bed dam-break problem.
At each time, strong descending-gradient candidates identify the leading downstream transition; local states immediately behind and ahead of that transition define its $50\%$ height crossing, which is linearly interpolated between the adjacent samples.
The speed is the least-squares slope of the front positions at $t=1$, $2$, and $3\,\second$.
Squares denote perfect slip, triangles denote no slip, and the dashed and dotted horizontal lines show the corresponding OpenFOAM values.
Under perfect slip, the reduced models overpredict the OpenFOAM speed by approximately $8\%$.
Under no slip, HSWME with $N=1$ gives the closest speed, whereas the modified models overpredict by approximately $9\%$ despite their closer resolved-interior profiles.}}
    \label{fig:ex-new-model-speed}
\end{figure}
\FloatBarrier

{ Figure~\ref{fig:ex-new-model-speed} demonstrates that vertical-profile fidelity and front-speed fidelity are distinct diagnostics.
Although the modified models better represent the resolved interior profile, the classical HSWME with $N=1$ gives the closest no-slip speed in this test.
The effective wall source therefore reduces excessive frictional damping but does not guarantee improved propagation speed.
This motivates comparing the retained-order behavior of the classical and modified moment systems.}

\FloatBarrier

{ To distinguish retained-order effects from the source modification, we repeat the no-slip HSWME and MHSWME calculations for every retained order $N=3,4,\ldots,15$ and compare the representative orders $N=5$, $10$, and $15$ at $t=3\,\second$.
All numerical settings are unchanged from the preceding no-slip comparison: $N_x=4000$, $\kappa=10^{10}\,\kg\per (\meter^2\cdot\second)$, and the OpenFOAM reference uses the $4000\times640$ mesh and pure-water cells ($\alpha_\star=1$).}

\begin{figure}[htbp]
    \centering
    \includegraphics[width=0.96\linewidth]{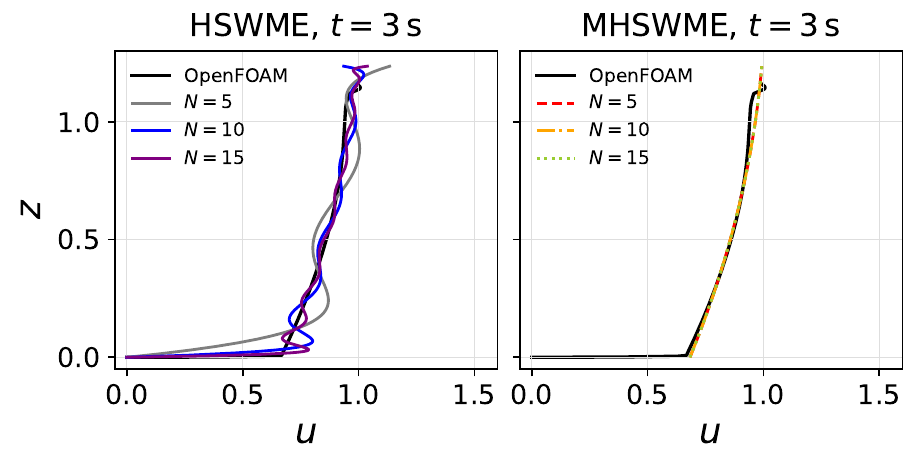}
    \caption{{ Retained-order comparison of the no-slip HSWME (left) and MHSWME (right): vertical profiles of horizontal velocity $u(z)$ at $x=55\,\meter$ and $t=3\,\second$ for $N=5$, $10$, and $15$.
The OpenFOAM reference uses the $4000\times640$ mesh and retains only pure-water cell centers ($\alpha_\star=1$), together with the exact no-slip wall trace $u(0)=0$.
As $N$ increases, the HSWME reconstruction narrows the modeled wall-to-interior transition but develops small oscillations.
The three MHSWME reconstructions are nearly indistinguishable and closely follow the resolved OpenFOAM interior.
They are not forced through the exact wall trace and do not resolve the thin physical wall layer.}}
    \label{fig:ex-new-hswme-order-profiles}
\end{figure}
\FloatBarrier

\begin{figure}[htbp]
    \centering
    \includegraphics[width=0.82\linewidth]{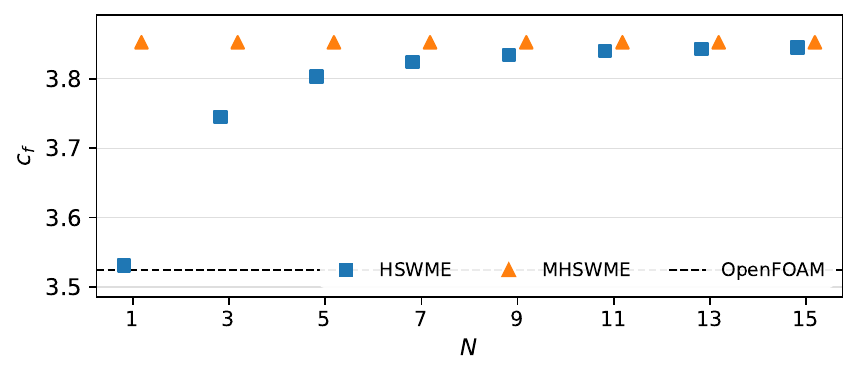}
    \caption{{ Effective right-moving front speed $c_f$ for the no-slip HSWME and MHSWME at retained orders $N=1,3,5,7,9,11,13$, and $15$.
For each order and time, the leading downstream transition is selected from strong descending-gradient candidates, its local $50\%$ height crossing is linearly interpolated between adjacent samples, and $c_f$ is the least-squares slope of the positions at $t=1$, $2$, and $3\,\second$.
The dashed line is the $4000\times640$ OpenFOAM value, $3.525\,\meter\per\second$.
The HSWME speed increases from $3.531\,\meter\per\second$ at $N=1$ to $3.845\,\meter\per\second$ at $N=15$, moving away from OpenFOAM over the tested orders.
The MHSWME speed remains essentially order-independent at $3.853\,\meter\per\second$, approximately $9.31\%$ above the OpenFOAM value.}}
    \label{fig:ex-new-hswme-order-speed}
\end{figure}
\FloatBarrier

{ Figure~\ref{fig:ex-new-hswme-order-profiles} shows two distinct retained-order behaviors at $t=3\,\second$.
Increasing the HSWME order sharpens the wall-to-interior transition over the displayed interval, but the global Legendre reconstruction develops small oscillations.
The MHSWME profiles are nearly independent of the retained order and overlap the OpenFOAM curve over most of the common resolved interior.
The mismatch at $z=0$ has a distinct meaning: the OpenFOAM curve includes the exact physical trace $u(0)=0$, whereas each plotted moment curve is the truncated Legendre reconstruction evaluated at $\zeta=0$; the MHSWME closure represents the unresolved wall influence without imposing that exact endpoint value on the reconstruction.}

{ Figure~\ref{fig:ex-new-hswme-order-speed} exhibits the corresponding contrast in front speed.
The MHSWME values are essentially unchanged with $N$, consistent with the nearly coincident vertical profiles in Figure~\ref{fig:ex-new-hswme-order-profiles}.
The HSWME speed instead increases with $N$ and approaches the MHSWME value over the tested orders, although its $N=1$ value is closest to OpenFOAM.
Thus, the sharper local HSWME vertical representation at higher order does not translate into a more accurate cumulative horizontal propagation speed in this test.
No profile convergence rate, integral profile error, or asymptotic speed convergence is claimed from these finite-order comparisons.}

{ The additional cost increases with the retained order.
In this one-dimensional implementation, the HSWME with order $N$ has $N+2$ state variables per cell.
With $N_x=4000$ and CFL $=0.7$, all tested orders require 741--747 time steps, while $N=5$, $10$, and $15$ take approximately $2.8$, $3.9$, and $5.3$ times the wall-clock time of $N=1$, respectively.
The state storage at $N=15$ is $17/3\approx5.7$ times that at $N=1$.
These ratios are implementation-dependent and are reported to indicate the additional cost rather than as universal performance measurements.}

\FloatBarrier
}

\begin{example}[3D radial water collapse]\label{ex2}
    \begin{revisionblock}
    We consider the 3D radial water collapse.
    The computational domain is $[0,100]^2\times[0,2]\,\meter^3$, where the $z$-direction is used only for the OpenFOAM simulations of the incompressible Navier--Stokes equations and is not needed for the shallow water models.
    We choose the initial water column
    \begin{equation}
    h(0,x,y) = \begin{cases}
        1.5 \,\meter, \quad \sqrt{(x-50)^2+(y-50)^2} \leq 15 \,\meter, \\
        1 \,\meter, \quad \mbox{ otherwise},
    \end{cases}
\end{equation}
with a constant zero initial velocity profile.
In OpenFOAM, the bottom and four lateral boundaries are no-slip walls, while the upper boundary is the atmosphere; $p_{\rm rgh}$ and $\alpha_{\rm w}$ use the wall and atmospheric conditions stated above.
The reduced models instead use zero-normal-gradient copies on all four lateral boundaries and represent the bottom through the source term, with $\kappa=10^{10}\,\kg\per(\meter^2\cdot\second)$ as a finite approximation to the no-slip limit.
Thus, the lateral boundary conditions are not identical between the benchmark and reduced problems; the comparisons below concern the interior collapse over $0\leq t\leq3\,\second$.
    \end{revisionblock}
\end{example}

{ We first assess vertical-resolution sensitivity of the OpenFOAM benchmark at fixed $N_x=N_y=400$ using $N_z=200$, $400$, and $800$.
Similar simulations of 3D radial water collapse have also been considered in \cite{Verbiest2025}.}

\begin{figure}[htbp]
    \centering
    \includegraphics[width=0.96\linewidth]{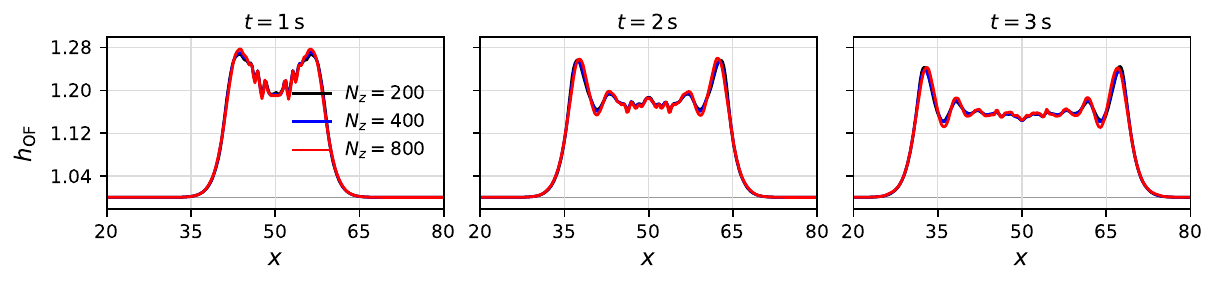}
    \caption{{ OpenFOAM vertical-resolution study of the water-height slice $h_{\rm OF}(t,x)$ at $y=65.5\,\meter$ for the 3D radial water collapse.
The panels show $t=1$, $2$, and $3\,\second$ from left to right; $N_x=N_y=400$ and the curves use $N_z=200$, $400$, and $800$.
The bulk front locations and large-scale profiles are similar across the meshes, while the local oscillatory free-surface structure retains visible vertical-resolution dependence.
The $N_z=800$ result is used as the finest available reference, without claiming pointwise convergence of the local oscillations.}}
    \label{fig:ex2-openfoam-hslice}
\end{figure}

\begin{figure}[htbp]
    \centering
    \includegraphics[width=0.96\linewidth]{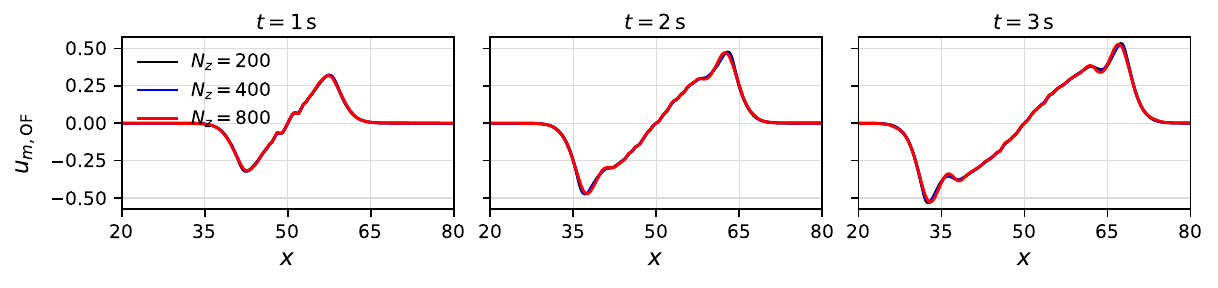}
    \caption{{ OpenFOAM vertical-resolution study of the depth-averaged horizontal velocity $u_{m,{\rm OF}}(t,x)$ at $y=65.5\,\meter$ for the 3D radial water collapse.
The panels show $t=1$, $2$, and $3\,\second$ from left to right; $N_x=N_y=400$ and the curves use $N_z=200$, $400$, and $800$.
The three curves are nearly coincident across the displayed interval, showing that this bulk velocity diagnostic is substantially less sensitive to the tested vertical refinement than the local water-height structure.}}
    \label{fig:ex2-openfoam-meanvelocity}
\end{figure}

{ Figures~\ref{fig:ex2-openfoam-hslice} and \ref{fig:ex2-openfoam-meanvelocity} distinguish the resolution sensitivity of local and bulk OpenFOAM diagnostics.
The local amplitudes and fine oscillatory structure near the moving free-surface transitions remain visibly mesh-dependent, so they are not interpreted as pointwise-converged features.
In contrast, the depth-averaged velocity is nearly unchanged over the tested meshes.
We therefore use $N_z=800$ in the remaining comparisons because it is the finest available calculation, not because every local interface feature has been proved converged.}

\begin{figure}[htbp]
    \centering
    \includegraphics[width=0.96\linewidth]{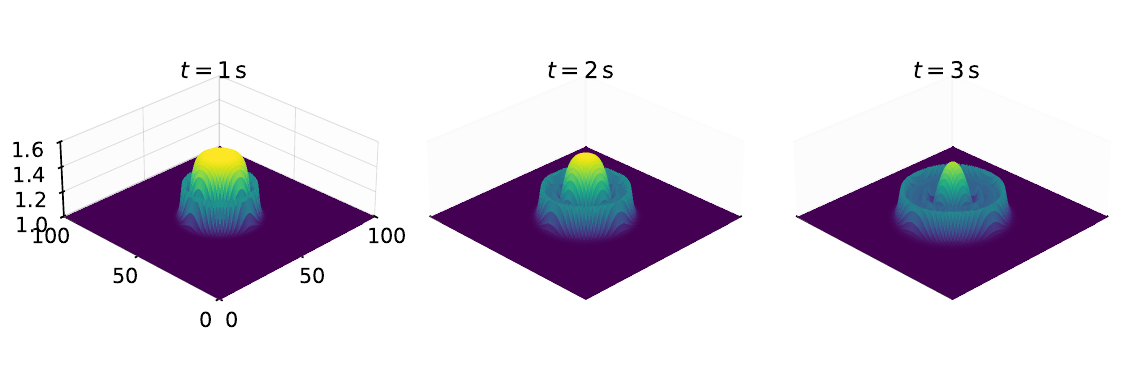}
    \caption{{ OpenFOAM water height (meters) over horizontal coordinates $(x,y)$ in meters for the 3D radial water collapse with zero initial velocity on the $400\times400\times800$ mesh at $t=1$, $2$, and $3\,\second$ from left to right.
The initially elevated circular column collapses into an approximately radially symmetric flow, producing an expanding annular front and a progressively lower central region.}}
    \label{fig:ex2-NS-surface}
\end{figure}

\begin{figure}[htbp]
    \centering
    \includegraphics[width=0.9\linewidth]{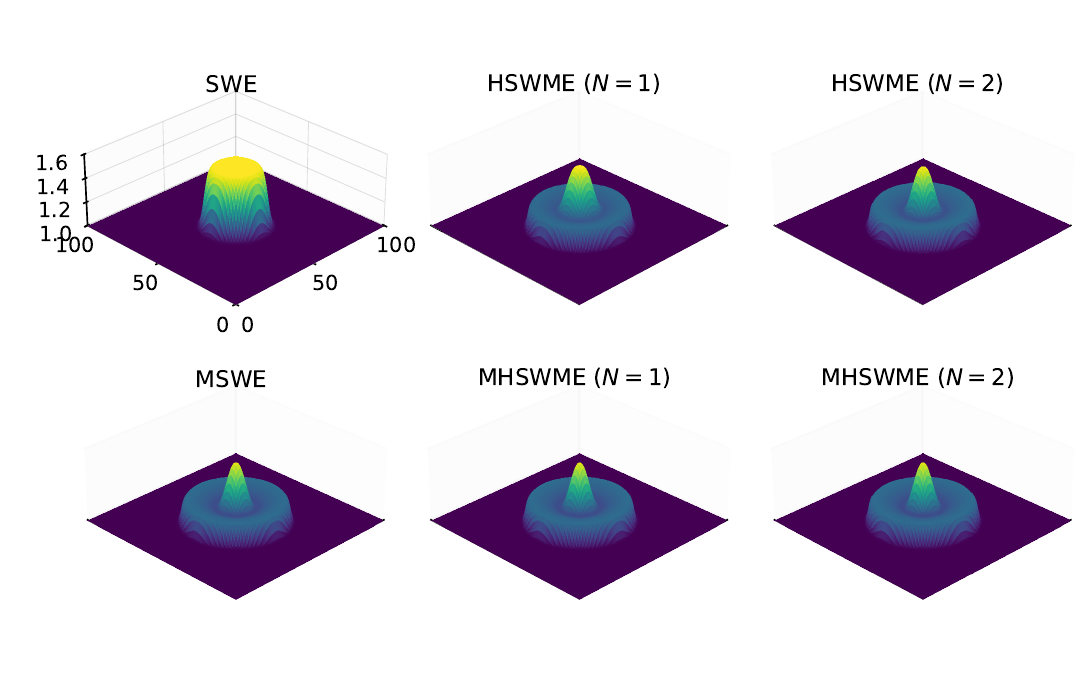}
    \caption{{ Water height at $t=3\,\second$ for the reduced models in the 3D radial water collapse: SWE, HSWME with $N=1,2$, MSWE, and MHSWME with $N=1,2$.
The SWE retains a pronounced central column and predicts much less collapse and outward spreading.
All moment-based and modified models produce a substantially lower center and an outward-propagating annular structure, with broadly similar large-scale surfaces at the displayed resolution.}}
    \label{fig:ex2-surface}
\end{figure}

\begin{figure}[htbp]
    \centering
    \includegraphics[width=0.9\linewidth]{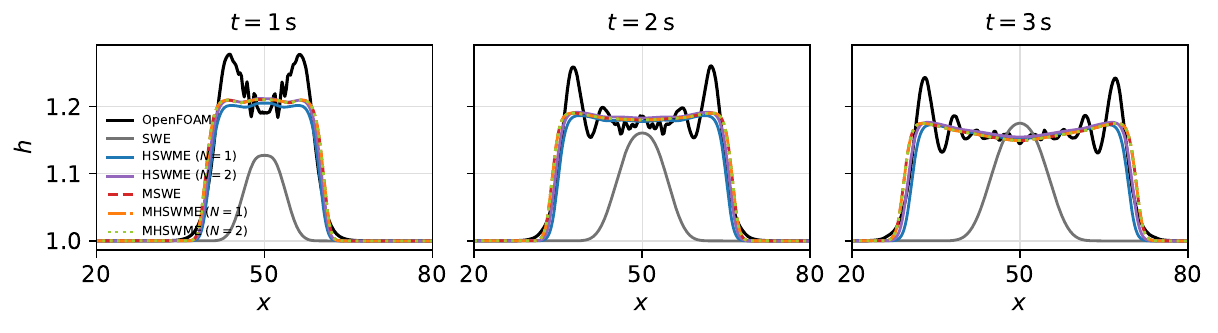}
    \caption{{ Water-height slices $h(t,x)$ at $y=65.5\,\meter$ for the 3D radial water collapse at $t=1$, $2$, and $3\,\second$ from left to right.
Curves compare the $400\times400\times800$ OpenFOAM reference, SWE, HSWME with $N=1,2$, MSWE, and MHSWME with $N=1,2$.
The SWE front propagates too slowly and retains a central hump, whereas the other reduced models capture the broad outward motion and central depression seen in OpenFOAM.
None of the reduced models reproduces the localized OpenFOAM oscillations near the moving fronts.}}
    \label{fig:ex2-height-y65}
\end{figure}

\begin{figure}[htbp]
    \centering
    \includegraphics[width=0.9\linewidth]{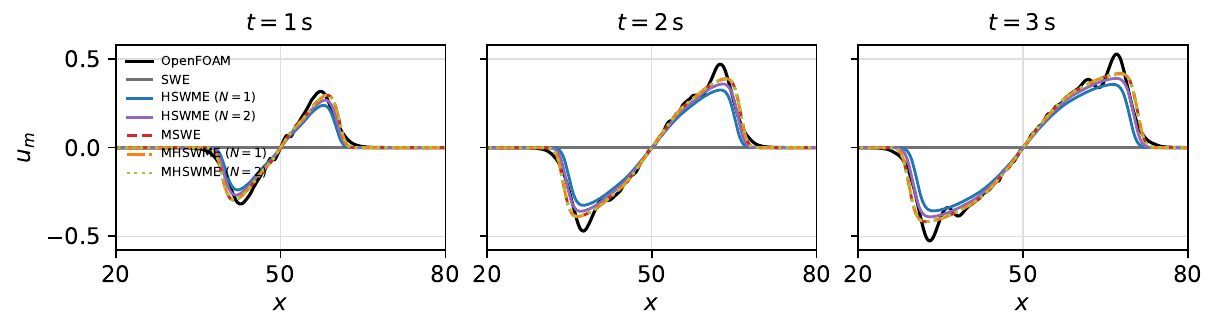}
    \caption{{ Depth-averaged horizontal velocity $u_m(t,x)$ at $y=65.5\,\meter$ for the 3D radial water collapse at $t=1$, $2$, and $3\,\second$ from left to right.
Curves compare the $400\times400\times800$ OpenFOAM reference, SWE, HSWME with $N=1,2$, MSWE, and MHSWME with $N=1,2$.
The approximately radial geometry produces an antisymmetric profile along the slice.
The SWE remains near zero, while the HSWME and modified models capture the sign, bulk magnitude, and outward propagation of the OpenFOAM velocity; differences among these reduced models are smaller than their common discrepancy from the sharp OpenFOAM extrema at the fronts.}}
    \label{fig:ex2-meanvelocity-y65}
\end{figure}

\begin{figure}[htbp]
    \centering
    \includegraphics[width=0.8\linewidth]{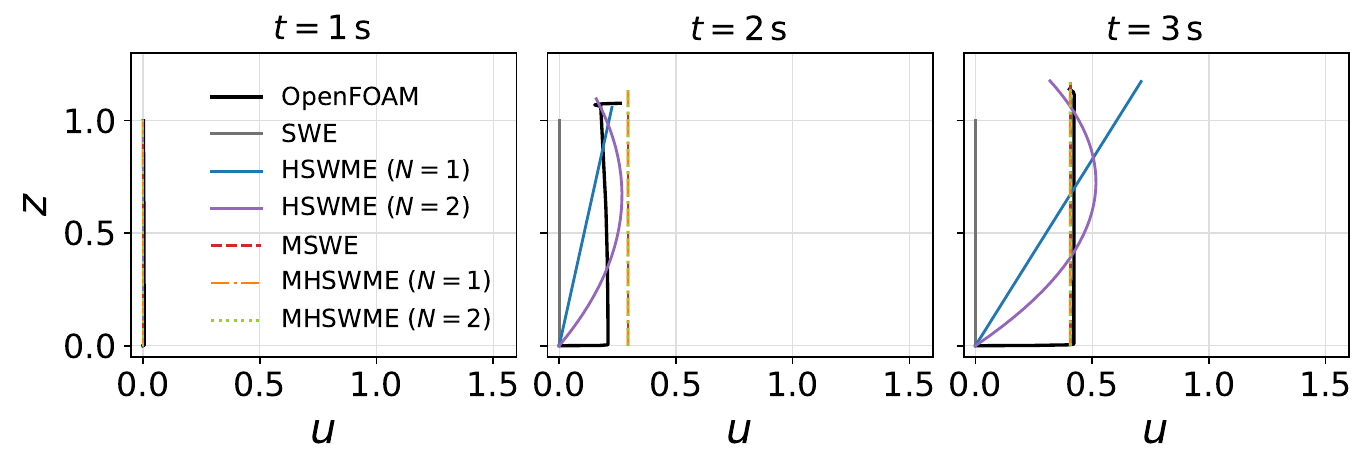}
    \caption{{ Vertical profiles of horizontal velocity $u(z)$ at $(x,y)=(65.5,65.5)\,\meter$ for the 3D radial water collapse at $t=1$, $2$, and $3\,\second$ from left to right.
OpenFOAM uses pure-water cells only ($\alpha_\star=1$); reduced profiles compare SWE, HSWME with $N=1,2$, MSWE, and MHSWME with $N=1,2$.
As the front reaches the sample point, OpenFOAM develops a thin wall-affected region and an approximately uniform resolved interior.
The low-order HSWME reconstructions pass through the zero wall endpoint and thereby distort the interior profile, while the modified models follow the resolved interior more closely.
No reduced curve is claimed to resolve the physical wall layer or the VOF interface.}}
    \label{fig:ex2-verticalvelocity-y65}
\end{figure}

\begin{figure}[htbp]
    \centering
    \includegraphics[width=0.8\linewidth]{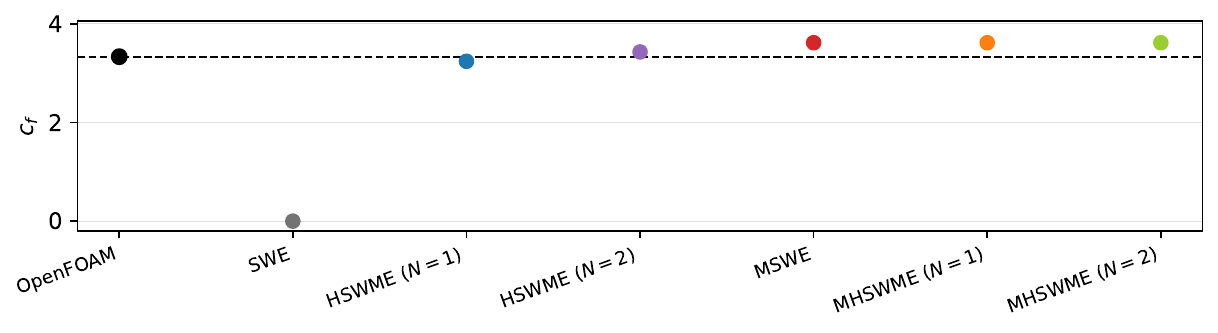}
    \caption{{ Effective outward front speed $c_f$ along the positive-$x$ centerline $y=50\,\meter$ for the 3D radial water collapse.
At each time, strong descending-gradient candidates identify the leading outward transition; local states behind and ahead define its $50\%$ height crossing, which is linearly interpolated between adjacent samples, and $c_f$ is the least-squares slope at $t=1$, $2$, and $3\,\second$.
The dashed line is the $400\times400\times800$ OpenFOAM value, $3.332\,\meter\per\second$.
The SWE is nearly arrested; the HSWME values for $N=1$ and $2$ bracket the OpenFOAM speed with approximately $2.8\%$ and $2.9\%$ absolute relative error, while the modified models overpredict by approximately $8.5\%$.
Thus, the modified closure improves the displayed interior profiles but not this front-speed diagnostic.}}
    \label{fig:ex2-propagation-speed}
\end{figure}

{ Figures~\ref{fig:ex2-NS-surface}--\ref{fig:ex2-propagation-speed} compare the large-scale surface evolution, horizontal slices, vertical structure, and effective front speed for Example~\ref{ex2}; the prescribed state at $t=0$ is omitted.
Because the problem is approximately symmetric in $x$ and $y$, only one set of equivalent horizontal slices is shown.}

{ The surface sequence provides the higher-fidelity benchmark for the large-scale collapse geometry.
Its nearly radial annular front motivates the one-dimensional slice comparisons below, while no claim of exact radial symmetry is made for the Cartesian discretization.}

{ The surface comparison shows that the classical SWE is dominated by the large wall-friction source and therefore remains close to the initial column.
Adding moments or using the effective wall source restores the broad collapse and outward annular motion seen in OpenFOAM.
Because these surfaces provide no water-height error norm, they support this qualitative separation but not a unique ranking among the moment-based and modified models.}

{ The slice at $y=65.5\,\meter$, just outside the initial column, resolves the timing difference more clearly than the surface views.
OpenFOAM has already passed the sampling line at $t=1\,\second$, whereas the HSWME and modified fronts are just arriving and the SWE remains strongly delayed.
At later times, all models except the SWE reproduce the broad central depression and outward motion, but their hydrostatic reduced dynamics do not reproduce the localized OpenFOAM structures near the fronts.}

{ The sign change reflects outward motion on the two sides of the nearly radial collapse, while the sharp OpenFOAM extrema coincide with the localized front structures in Figure~\ref{fig:ex2-height-y65}.
The moment-based and modified models give broadly comparable bulk velocities away from those extrema.
Their propagation accuracy must nevertheless be assessed separately because similar mean-velocity shapes do not imply identical front speeds.}

{ The small velocity at $t=1\,\second$ is consistent with the front just reaching the sampling point.
At later times, no slip creates a strong wall-to-interior gradient, while the resolved pure-water interior remains comparatively uniform.
Enforcing the zero endpoint through a low-order global HSWME polynomial then perturbs the entire reconstructed profile; the modified source represents the unresolved wall influence through \eqref{eq:modified-gamma} without imposing that endpoint on the polynomial.
When shown, the physical zero wall value is a separate visualization datum, not an additional correction to the modified reconstruction, and the VOF interface is excluded from the pointwise profile assessment.}

\FloatBarrier

{ This speed comparison confirms that the effective source substantially avoids the near-arrest produced by the classical SWE source, but it also exposes a diagnostic tradeoff.
The two HSWME speeds are closer to OpenFOAM even though their low-order vertical reconstructions are less representative of the resolved interior.
Consequently, neither the surface shapes nor the vertical profiles alone determine propagation accuracy.}

{ Taken together, the diagnostics in Example~\ref{ex2} show that the classical SWE source excessively damps the collapse, while adding moments recovers the broad height and mean-velocity evolution but leaves the low-order reconstructed interior profile sensitive to the zero wall endpoint.
The modified source follows the resolved interior profile more closely and avoids near-arrest, yet its centerline front speed is farther from OpenFOAM than the two HSWME values.
The improvement is therefore diagnostic-specific rather than a uniform accuracy ordering.}

\begin{example}[{ 3D water collapse with linear initial velocity}]\label{ex:3d-nonzero}
    \begin{revisionblock}
    We use the same geometry, water height, and no-slip bottom condition as in Example~\ref{ex2}, but deliberately prescribe the linear initial velocity profile
    \begin{equation}\label{init:ex3}
        u(0,x,y,z)=v(0,x,y,z)=Sz,
    \end{equation}
with $S=0.25 \,\second^{-1}$.
This profile satisfies the bottom value $u=v=0$ but is not compatible with the stress-free free-surface condition used in the reduced derivation, because
\begin{equation}\label{eq:ex3-free-surface-incompatibility}
    \partial_z u=\partial_z v=S\neq0
    \qquad \mbox{at }z=h(0,x,y).
\end{equation}
It is therefore used as a deliberately sheared nonequilibrium initial state, not as data satisfying every compatibility condition of the free-surface initial-boundary-value problem.

In OpenFOAM, $x=0$ and $y=0$ maintain a $1\,\meter$ water layer with $u=v=Sz$, while $x=100\,\meter$ and $y=100\,\meter$ use zero-normal-gradient outflow; the bottom is no slip and the upper boundary is the atmosphere.
The same incompatible shear profile is therefore continuously introduced at the two inflows rather than appearing only at $t=0$.
For the reduced models, the exact Legendre projection at a local depth $h$ is
\begin{equation}\label{eq:ex3-initial-moments}
    u_m=v_m=\frac{Sh}{2},
    \qquad
    \ap_1=\be_1=-\frac{Sh}{2},
    \qquad
    \ap_i=\be_i=0\quad (i\geq2),
\end{equation}
and this projected state with $h=1\,\meter$ is fixed at $x=0$ and $y=0$; the opposite boundaries use zero-normal-gradient outflow.
We again use $\kappa=10^{10}\,\kg\per(\meter^2\cdot\second)$ as a finite approximation to the no-slip limit.
Consequently, the results at $t=1$, $2$, and $3\,\second$ include adjustment from both the free-surface-incompatible startup and the maintained sheared inflow, and the conclusions are restricted to this forced test.
    \end{revisionblock}
\end{example}

\begin{figure}[htbp]
    \centering
    \includegraphics[width=0.96\linewidth]{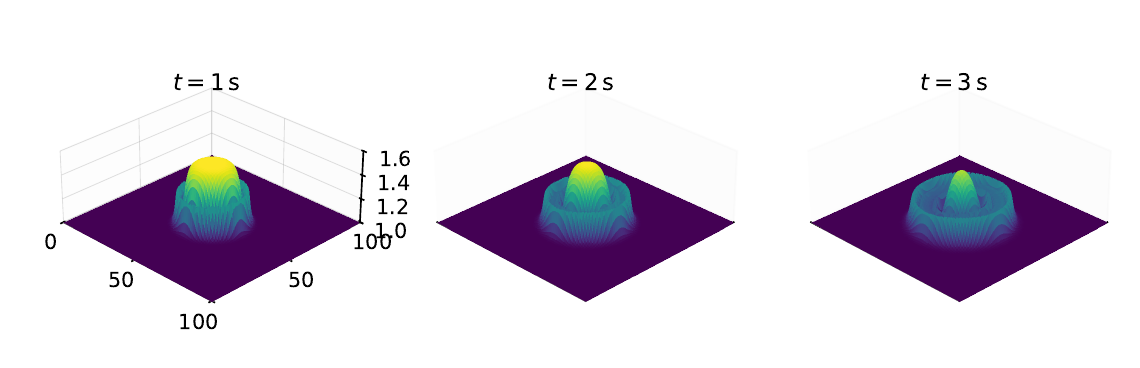}
    \caption{{ OpenFOAM water height (meters) over horizontal coordinates $(x,y)$ in meters for the 3D collapse with the deliberately sheared initial velocity $u(0,x,y,z)=v(0,x,y,z)=Sz$, $S=0.25\,\second^{-1}$, on the $400\times400\times800$ mesh at $t=1$, $2$, and $3\,\second$ from left to right.
The column collapses and spreads outward, but the positive-$x$ and positive-$y$ initial velocities together with the maintained sheared inflows break radial symmetry, producing a higher and faster-moving front near the $x=y$ direction than in the opposite directions.
The sequence includes the adjustment caused by the free-surface incompatibility in~\eqref{eq:ex3-free-surface-incompatibility}.}}
    \label{fig:ex3-NS-surface}
\end{figure}

\begin{figure}[htbp]
    \centering
    \includegraphics[width=0.96\linewidth]{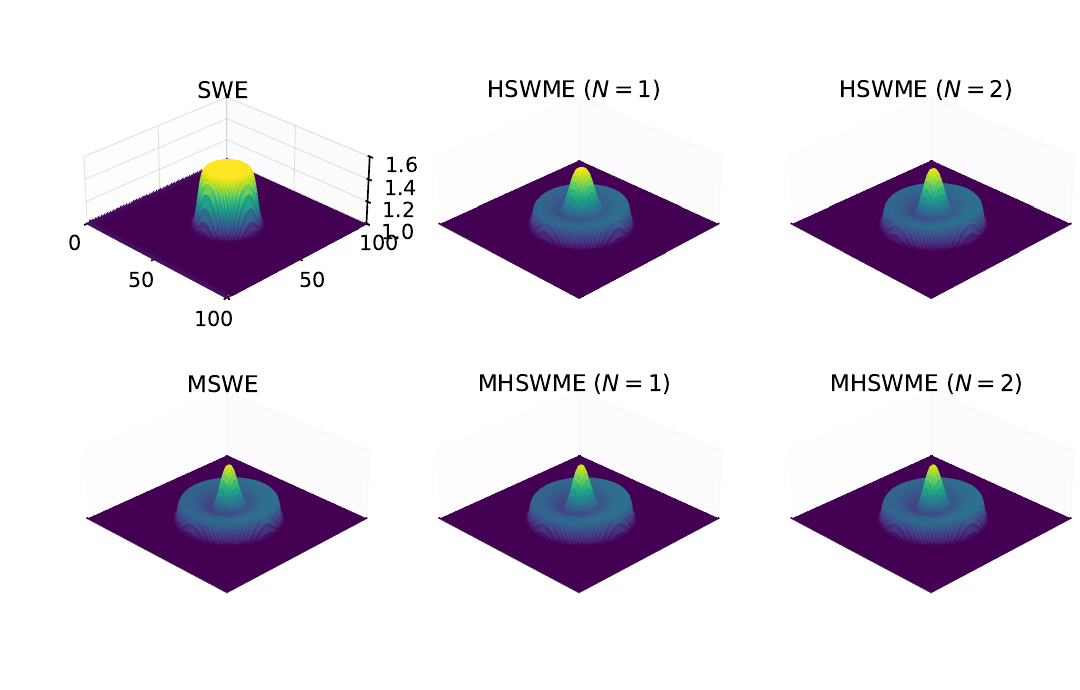}
    \caption{{ Water height at $t=3\,\second$ for the reduced models in the 3D collapse with linear initial velocity: SWE, HSWME with $N=1,2$, MSWE, and MHSWME with $N=1,2$.
The SWE retains an overly tall central column and strongly underpredicts the collapse.
The other models produce similar large-scale outward spreading and reproduce the directional asymmetry induced by the initial and inflow velocities.
These surfaces alone do not establish a unique accuracy ranking because no water-height error norm is shown.}}
    \label{fig:ex3-surface}
\end{figure}

\begin{figure}[htbp]
    \centering
    \includegraphics[width=0.96\linewidth]{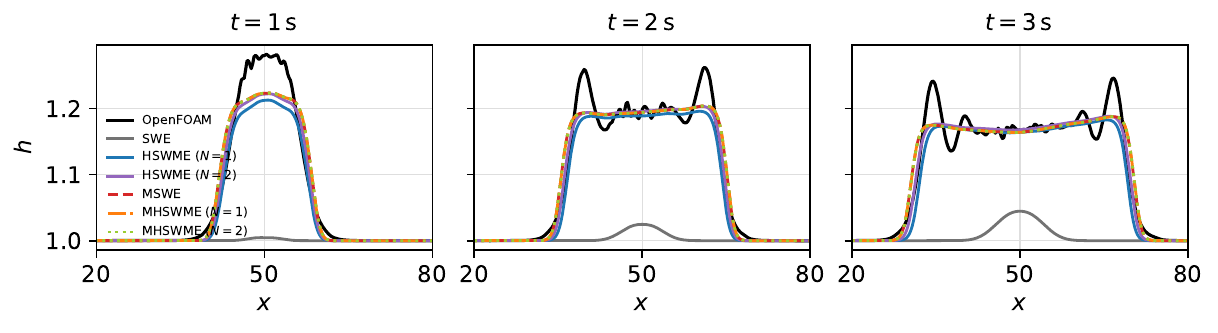}
    \caption{{ Water-height slices $h(t,x)$ at $y=67\,\meter$ for the 3D collapse with linear initial velocity at $t=1$, $2$, and $3\,\second$ from left to right.
Curves compare OpenFOAM, SWE, HSWME with $N=1,2$, MSWE, and MHSWME with $N=1,2$.
The SWE front remains substantially delayed and the slice retains a weak central hump.
The other reduced models capture the broad OpenFOAM front motion and the stronger positive-going branch produced by the imposed velocity, but not the localized OpenFOAM oscillations at the fronts.}}
    \label{fig:ex3-height-y67}
\end{figure}

\begin{figure}[htbp]
    \centering
    \includegraphics[width=0.96\linewidth]{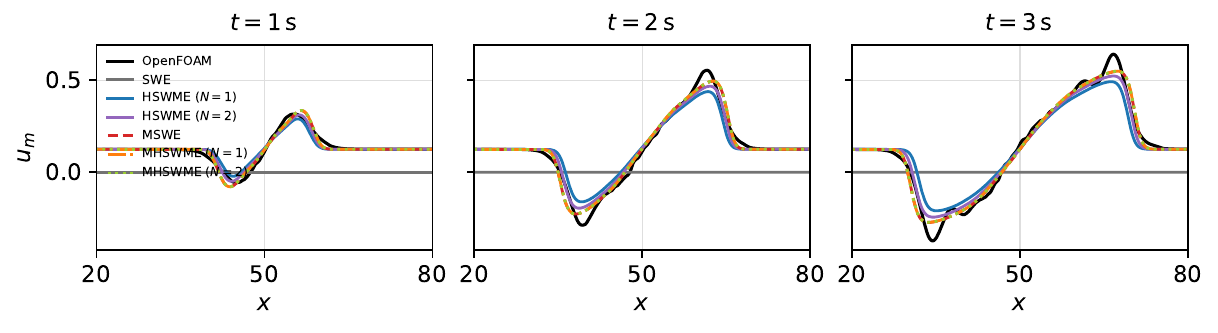}
    \caption{{ Depth-averaged horizontal velocity $u_m(t,x)$ at $y=67\,\meter$ for the 3D collapse with linear initial velocity at $t=1$, $2$, and $3\,\second$ from left to right.
Curves compare OpenFOAM, SWE, HSWME with $N=1,2$, MSWE, and MHSWME with $N=1,2$.
The imposed profile $u=v=Sz$, with $S=0.25\,\second^{-1}$, and the inflow make the positive-going front stronger than the negative-going front, so the slice is not point-symmetric.
The SWE remains nearly stationary, whereas the HSWME and modified models capture the bulk asymmetric OpenFOAM velocity; the sharpest OpenFOAM extrema near the fronts remain under-resolved.}}
    \label{fig:ex3-meanvelocity-y67}
\end{figure}

\begin{figure}[htbp]
    \centering
    \includegraphics[width=0.8\linewidth]{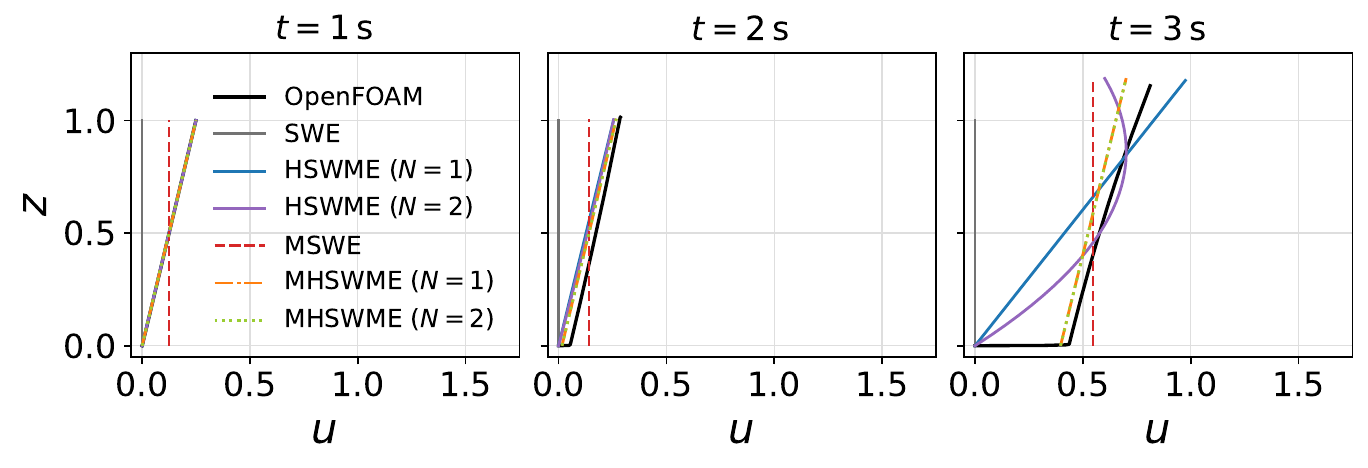}
    \caption{{ Vertical profiles of horizontal velocity $u(z)$ at $(x,y)=(67,67)\,\meter$ for the 3D collapse with linear initial velocity at $t=1$, $2$, and $3\,\second$ from left to right.
OpenFOAM uses pure-water cells only ($\alpha_\star=1$); reduced profiles compare SWE, HSWME with $N=1,2$, MSWE, and MHSWME with $N=1,2$.
At $t=1$ and $2\,\second$, the wall-to-interior velocity difference is small and the low-order HSWME profiles remain close to the nearly linear OpenFOAM interior.
By $t=3\,\second$, that difference is larger, and enforcing the zero endpoint through a low-order global HSWME polynomial visibly distorts the interior.
The modified models remain closer to the resolved OpenFOAM interior, while neither model class resolves the thin physical wall layer.}}
    \label{fig:ex3-verticalvelocity-y67}
\end{figure}

\begin{figure}[htbp]
    \centering
    \includegraphics[width=0.8\linewidth]{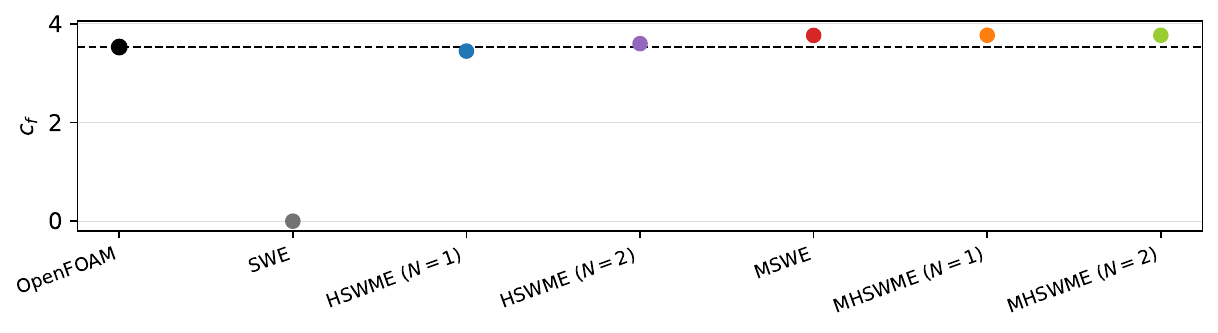}
    \caption{{ Effective positive-$x$ front speed $c_f$ along $y=50\,\meter$ for the 3D collapse with linear initial velocity.
At each time, strong descending-gradient candidates identify the leading positive-going transition; local states behind and ahead define its $50\%$ height crossing, which is linearly interpolated between adjacent samples, and $c_f$ is the least-squares slope at $t=1$, $2$, and $3\,\second$.
The dashed line is the $400\times400\times800$ OpenFOAM value, $3.528\,\meter\per\second$.
The SWE is nearly stationary; the HSWME values for $N=1$ and $2$ have approximately $2.3\%$ and $2.0\%$ absolute relative error, while the modified models overpredict by approximately $6.7$--$6.8\%$.
Thus, closer resolved-interior profile agreement does not imply the closest front speed.}}
    \label{fig:ex3-propagation-speed}
\end{figure}

{ Figures~\ref{fig:ex3-NS-surface}--\ref{fig:ex3-propagation-speed} compare the large-scale surface evolution, asymmetric horizontal slices, vertical structure, and positive-$x$ front speed for Example~\ref{ex:3d-nonzero}; the prescribed state at $t=0$ is omitted.
These diagnostics include the subsequent adjustment of the incompatible initial shear and the continuing influence of the two prescribed inflows.}

{ Relative to the zero-initial-velocity case, the prescribed positive horizontal motion biases both the initial discharge and the subsequent boundary forcing.
The resulting directional asymmetry provides a distinct test of whether the reduced models reproduce the bulk transport rather than only the radial-collapse geometry.
It should not, however, be interpreted as the response to free-surface-compatible equilibrium data.}

{ As in the radial case, the classical SWE remains dominated by excessive wall-source damping, whereas the moment-based and modified models recover the broad collapse.
Their asymmetric fronts also show that the imposed positive-direction transport is retained by the reduced dynamics.
Because the surfaces are visually similar among these models and no norm is computed, the slice and profile diagnostics below are needed for further discrimination.}

{ Because $y=67\,\meter$ lies farther from the initial column than the slice used in Example~\ref{ex2}, this comparison emphasizes arrival time and directional transport.
The positive-going front is stronger than the negative-going front, consistently with the imposed velocity and inflow.
The remaining reduced models reproduce that bulk asymmetry, while the SWE delay persists despite the nonzero initial motion.}

{ The loss of point symmetry separates the effect of the imposed directional transport from the nearly radial geometry of Example~\ref{ex2}.
Away from the sharp front extrema, the moment-based and modified models recover similar bulk velocity levels and preserve the nonzero background motion generated by the initial and inflow conditions.
The SWE instead remains close to zero because its wall source overwhelms that transport.}

{ This time sequence shows that the severity of the endpoint constraint depends on the wall-to-interior velocity difference, not only on the polynomial order.
When that difference is small, the classical low-order profiles remain visually close to the resolved interior; as it grows, the endpoint condition distorts the global reconstruction.
The modified profiles then remain visually closer over the common pure-water interval, although no profile norm or pointwise wall-layer resolution is claimed.}

\FloatBarrier

{ The front-speed ordering again differs from the vertical-profile ordering.
The HSWME speeds are closer to OpenFOAM even though the modified profiles remain visually closer to the resolved interior at $t=3\,\second$.
The effective source therefore changes the damping and vertical reconstruction in a way that is beneficial for some diagnostics but not uniformly for propagation.}

{ Taken together, Example~\ref{ex:3d-nonzero} shows that the classical SWE source suppresses the imposed transport, whereas the moment-based and modified models capture the broad asymmetric height and mean-velocity evolution.
The classical low-order vertical profiles remain visually close to OpenFOAM while the wall-to-interior difference is small but become globally distorted as that difference grows.
The modified profiles then remain closer over the resolved interior, while HSWME with $N=2$ gives the closest front speed.
These conclusions are diagnostic-specific and restricted to this deliberately forced, free-surface-incompatible test, the stated boundary conditions, and $\kappa=10^{10}\,\kg\per(\meter^2\cdot\second)$; they do not establish accuracy for a fully compatible free-surface initial-boundary-value problem.}

\section{Conclusions}
\label{sec:conclusion}

{ We developed an effective wall-traction closure for the SWE and the specified HSWME principal part by replacing the reconstructed-endpoint wall source with a bottom-to-mean relation derived in the distinguished regular-friction regime.
The closure represents the influence of unresolved near-wall dynamics without resolving the physical no-slip boundary layer or imposing the exact wall trace on a finite-order Legendre reconstruction.
Because only the source term is modified, the homogeneous principal matrices and their two-dimensional hyperbolicity classification remain unchanged.

Example~\ref{ex:2d-slip-noslip} provides a limiting control: under perfect slip, each classical model coincides with its modified counterpart, confirming that the modification is inactive when the wall traction vanishes.
Under the large-friction approximation to no slip, the classical low-order wall source produces excessive damping and can distort the resolved-interior profile, whereas the modified source recovers more of the broad water-height, depth-averaged-velocity, and interior-profile behavior displayed by OpenFOAM.
As the HSWME order increases, its vertical-profile representation improves but develops small oscillations, while its front speed moves away from the OpenFOAM value over the tested orders.
The MHSWME profiles and front speeds remain nearly independent of the retained order, indicating that the principal benefit of the modified closure occurs in the low-order regime.
The three-dimensional examples show the same diagnostic distinction: the modified source better represents the resolved interior, while the HSWME front speeds are closer to OpenFOAM.
Consequently, neither the source modification nor increasing the retained order produces a uniform improvement across water height, mean velocity, vertical profile, and front speed.

These results remain model-to-benchmark comparisons because OpenFOAM and the reduced models solve different continuum problems, and no uniform model-accuracy ordering or pointwise profile convergence is claimed.
The $\mathcal{O}(\varepsilon^2)$ estimate applies only to the wall-source reconstruction under the smooth regular-friction assumptions, whereas the large-$\kappa$ calculations use the rational coefficient~\eqref{eq:modified-gamma} as an empirically assessed wall-model continuation.
The linear-shear example is additionally restricted to its deliberately forced, free-surface-incompatible initial and boundary data.
Accordingly, extension to other friction regimes, moment regularizations, or flow configurations requires separate analysis and validation.}

\section*{Acknowledgements}

The authors acknowledge support from  ONR grant N00014-24-1-2242, NSF grant DMS-2309655, NSF grant DMS-2618114, AFOSR grants FA9550-24-1-0254, and DOE grant DE-SC0023164.

\appendix
\renewcommand{\thesection}{Appendix~\Alph{section}}
\renewcommand{\thesubsection}{\Alph{section}.\arabic{subsection}}
\section{{ Other wall-friction regimes}}
\label{app:other-friction-regimes}

{
Section~\ref{sec:flat-asymptotic-correction} uses the distinguished scaling $\gamma=\varepsilon\gamma_0$, for which the wall source and the inertial shallow-water terms are both $\mathcal{O}(1)$.
For completeness, this appendix records the other three formal friction scalings while retaining $\iRe=\mathcal{O}(1)$.
The same rational expression
\begin{equation}\label{eq:appendix-common-gamma-bar}
    \bar\gamma
    =\fl{\gamma}{\varepsilon\left(1+\fl{h\gamma}{3\iRe}\right)},
\end{equation}
appears in each reduced local vertical problem, but its order, the velocity scale, and the limiting evolution equation change with the friction regime.
}

\subsection{{ Very weak friction: $\gamma=\varepsilon^2\gamma_2$}}
\label{app:very-weak-friction}

{
Assume
\begin{equation}\label{eq:appendix-very-weak-scaling}
    \gamma=\varepsilon^2\gamma_2,
    \qquad \gamma_2>0,
    \qquad \gamma_2=\mathcal{O}(1).
\end{equation}
The leading field $\boldsymbol{U}_0$ and the first coefficient $\boldsymbol{u}_1$ are both depth-uniform.
The leading equations are therefore the frictionless shallow-water equations,
\begin{equation}\label{eq:appendix-very-weak-leading}
\begin{aligned}
    \fl{\pa h}{\pa t}+\nabla_H\!\cdot(h\boldsymbol{U}_0)&=0,\\
    \fl{\pa\boldsymbol{U}_0}{\pa t}
    +(\boldsymbol{U}_0\!\cdot\nabla_H)\boldsymbol{U}_0
    +G\nabla_H(h+h_b)&=\boldsymbol{0}.
\end{aligned}
\end{equation}
The first wall-dependent curvature occurs in the $\varepsilon^2$ coefficient.
For the reduced local problem with zero free-surface shear through that order,
\begin{equation}\label{eq:appendix-very-weak-profile}
    \boldsymbol{u}_H^\varepsilon(\xi)
    =\left[1+\fl{\gamma}{\iRe}
    \left(\xi-\fl{\xi^2}{2h}\right)\right]
    \boldsymbol{u}_b^\varepsilon
    +\mathcal{O}(\varepsilon^3),
\end{equation}
and
\begin{equation}\label{eq:appendix-very-weak-mean}
    \boldsymbol{u}_m^\varepsilon
    =\left(1+\fl{h\gamma}{3\iRe}\right)
    \boldsymbol{u}_b^\varepsilon
    +\mathcal{O}(\varepsilon^3).
\end{equation}
In this regime
\begin{equation}\label{eq:appendix-very-weak-gamma-bar}
    \bar\gamma
    =\fl{\varepsilon\gamma_2}
    {1+\fl{\varepsilon^2h\gamma_2}{3\iRe}}
    =\varepsilon\gamma_2+\mathcal{O}(\varepsilon^3).
\end{equation}
Thus the rational denominator first changes the wall source at $\mathcal{O}(\varepsilon^3)$.
The full second-order Navier--Stokes profile also contains the order-$\varepsilon^2$ free-surface tangential-traction contribution; consequently, \eqref{eq:appendix-very-weak-profile} represents the wall-induced part of the second-order correction unless that additional contribution is retained or shown to vanish.
}

\subsection{{ Order-one friction: $\gamma=\gamma_*=\mathcal{O}(1)$}}
\label{app:order-one-friction}

{
Let $\gamma=\gamma_*>0$ be independent of $\varepsilon$.
The $\mathcal{O}(\varepsilon^{-1})$ vertical-diffusion equation, the free-surface zero-shear condition, and the bottom Navier condition imply
\begin{equation}\label{eq:appendix-order-one-U0}
    \boldsymbol{U}_0=\boldsymbol{0}.
\end{equation}
Hence the nontrivial horizontal velocity is $\mathcal{O}(\varepsilon)$ rather than $\mathcal{O}(1)$.
Introduce the slow time $\tau=\varepsilon t$ and write
\begin{equation}\label{eq:appendix-order-one-slow-scaling}
    \boldsymbol{u}_H^\varepsilon
    =\varepsilon\boldsymbol{V}_0+\mathcal{O}(\varepsilon^2).
\end{equation}
With $\boldsymbol{F}_*=G\nabla_H(h+h_b)$, the leading local problem gives
\begin{equation}\label{eq:appendix-order-one-profile}
    \boldsymbol{V}_0(\xi)
    =-\left[
    \fl{h}{\gamma_*}
    +\fl{1}{\iRe}\left(h\xi-\fl{\xi^2}{2}\right)
    \right]\boldsymbol{F}_*.
\end{equation}
Its depth average is
\begin{equation}\label{eq:appendix-order-one-mean}
    \boldsymbol{V}_m
    =-G\left(\fl{h}{\gamma_*}+\fl{h^2}{3\iRe}\right)
    \nabla_H(h+h_b),
\end{equation}
and the wall-to-mean relation is
\begin{equation}\label{eq:appendix-order-one-wall-mean}
    \boldsymbol{V}_m
    =\left(1+\fl{h\gamma_*}{3\iRe}\right)\boldsymbol{V}_b.
\end{equation}
The effective coefficient is $\bar\gamma=\mathcal{O}(\varepsilon^{-1})$, so an $\mathcal{O}(\varepsilon^2)$ wall-trace error is amplified to an $\mathcal{O}(\varepsilon)$ source error.
The regular-friction remainder estimate therefore does not carry over.
The leading slow evolution is the lubrication equation
\begin{equation}\label{eq:appendix-order-one-lubrication}
    \fl{\pa h}{\pa\tau}
    =\nabla_H\!\cdot\left[
    G\left(\fl{h^2}{\gamma_*}+\fl{h^3}{3\iRe}\right)
    \nabla_H(h+h_b)
    \right].
\end{equation}
}

\subsection{{ Large friction: $\gamma=\Gamma_{-1}/\varepsilon$}}
\label{app:large-friction}

{
Finally, assume
\begin{equation}\label{eq:appendix-large-friction-scaling}
    \gamma=\fl{\Gamma_{-1}}{\varepsilon},
    \qquad \Gamma_{-1}>0,
    \qquad \Gamma_{-1}=\mathcal{O}(1).
\end{equation}
Again $\boldsymbol{U}_0=\boldsymbol{0}$ and the post-relaxation velocity is $\mathcal{O}(\varepsilon)$.
With the slow scaling
\begin{equation}\label{eq:appendix-large-friction-velocity}
    \boldsymbol{u}_H^\varepsilon
    =\varepsilon\boldsymbol{V}_0
    +\varepsilon^2\boldsymbol{V}_1
    +\mathcal{O}(\varepsilon^3),
\end{equation}
the leading bottom condition is no slip, $\boldsymbol{V}_0(0)=\boldsymbol{0}$.
The leading profile and mean are
\begin{equation}\label{eq:appendix-large-friction-profile}
\begin{aligned}
    \boldsymbol{V}_0(\xi)
    &=-\fl{G}{\iRe}
    \left(h\xi-\fl{\xi^2}{2}\right)
    \nabla_H(h+h_b),\\
    \boldsymbol{V}_{0,m}
    &=-\fl{Gh^2}{3\iRe}\nabla_H(h+h_b).
\end{aligned}
\end{equation}
The bottom trace is $\boldsymbol{u}_b^\varepsilon=\mathcal{O}(\varepsilon^2)$, while the depth average is $\boldsymbol{u}_m^\varepsilon=\mathcal{O}(\varepsilon)$.
The effective coefficient has the expansion
\begin{equation}\label{eq:appendix-large-friction-gamma-bar}
    \bar\gamma
    =\fl{3\iRe}{\varepsilon h}
    -\fl{9(\iRe)^2}{h^2\Gamma_{-1}}
    +\mathcal{O}(\varepsilon),
\end{equation}
and the leading slow equation is
\begin{equation}\label{eq:appendix-large-friction-lubrication}
    \fl{\pa h}{\pa\tau}
    =\nabla_H\!\cdot\left[
    \fl{Gh^3}{3\iRe}\nabla_H(h+h_b)
    \right].
\end{equation}
Thus the saturated value in \eqref{eq:modified-gamma-large-friction-limit} has a leading no-slip Stokes interpretation when $\iRe=\mathcal{O}(1)$, but the inertial shallow-water model is not the leading asymptotic system in this regime.
}

\begin{remark}[{ Strong-friction asymptotics versus the numerical no-slip setting}]
{ The quantity $\boldsymbol{u}_m$ has the same meaning throughout this work: it is the depth-averaged horizontal velocity.
The estimate $\boldsymbol{u}_m=\mathcal{O}(\varepsilon)$ in Sections~\ref{app:order-one-friction}--\ref{app:large-friction} follows from the joint assumption $\iRe=\mathcal{O}(1)$ and the corresponding strong wall scaling; it does not follow from a large dimensional value of $\kappa$ alone.
In the dam-break calculations, $\kappa=10^{10}$ is used to make the slip length $\ell_s=\mu/\kappa$ negligible and hence to approximate the pointwise no-slip condition.
For water on the advective dam-break time scale, however, $\iRe/\varepsilon\ll1$, so vertical diffusion acts initially only through a thin bottom layer and the inertia-dominated bulk depth average may remain $\mathcal{O}(1)$.
Thus the numerical large-$\kappa$ setting is not the $\iRe=\mathcal{O}(1)$ lubrication limit above.
The rational source is used there as a wall-model continuation, and no uniform regular-friction remainder estimate is claimed.}
\end{remark}

\section{{ Proof of Theorem \ref{thm:HSWME}}}
\label{app:proof-theorem1}

We first recall the standard definition of hyperbolicity in 2D:
\begin{definition}[2D hyperbolicity]\label{def:2D-hyperbolicity}
Consider the HSWME \eqref{sys:HSWME-vect}.
For $\theta\in[0,2\pi)$, define the matrix
\begin{equation}\label{eq:An-def}
A_n(\bbU,\theta):=\cos\theta\,A^{N}_{H}(\bbU)+\sin\theta\,B^{N}_{H}(\bbU).
\end{equation}
The system \eqref{sys:HSWME-vect} is {(strongly) hyperbolic in two dimensions at $\bbU$} if, for any $\theta\in[0,2\pi)$, the matrix $A_n(\bbU,\theta)$ is real diagonalizable.
It is {weakly hyperbolic in two dimensions at $\bbU$} if, for any $\theta\in[0,2\pi)$, the matrix $A_n(\bbU,\theta)$ has only real eigenvalues.
\end{definition}

To analyze the hyperbolicity in 2D, we use the rotational invariance of the HSWME in \cite{Bauerle2025}.
We emphasize that the real diagonalizability of $A_{H}^{N}(\bbU)$ at a fixed state is only a property of the one-dimensional coefficient matrix in the $x$-direction.
It is not, by itself, equivalent to the strong hyperbolicity of the two-dimensional system.
According to Definition \ref{def:2D-hyperbolicity}, strong 2D hyperbolicity requires the directional matrix $A_{n}(\bbU,\theta)$ to be real diagonalizable for every direction $\theta\in[0,2\pi)$.
The rotational invariance is used below to connect this directional condition to the diagonalizability of $A_{H}^{N}$ evaluated at all rotated states $T(\theta)\bbU$.

\begin{theorem}[Rotational invariance, Theorem 3.2 in \cite{Bauerle2025}]\label{thm:rotation-invariance-HSWME}
The HSWME \eqref{sys:HSWME-vect} satisfies the rotational invariance:
    \begin{equation}\label{eq:rotation-invariance-HSWME}
        \cos\theta \, A^{N}_{H}(\bbU) + \sin\theta  \, B^{N}_{H}(\bbU) = T^{-1} A^{N}_{H}(T \bbU) T.
    \end{equation}
    Here $T=T(\theta)$ is the rotation matrix given by
    \begin{equation}\label{eq:rotation-matrix-T}
        T(\theta) = \textrm{diag}(1, T_2(\theta), T_2(\theta), \cdots, T_2(\theta))\in\mathbb{R}^{(2N+3)\times(2N+3)},
    \end{equation}
    with
    \begin{equation}\label{eq:rotation-matrix-size-2-2}
    T_2(\theta) = 
    \begin{pmatrix}
    \cos\theta & \sin\theta \\
    -\sin\theta & \cos\theta 
    \end{pmatrix} 
    \in\mathbb{R}^{2\times 2}.
    \end{equation}    
\end{theorem}

\begin{proposition}[Reduction to $A_H^N$]\label{prop:reduction-AH}
For any fixed state $\bbU$ and any $\theta\in[0,2\pi)$, the matrices $A_n(\bbU,\theta)$ and $A_H^N(T(\theta)\bbU)$ are similar.
In particular,
\begin{equation}\label{eq:equiv-diag}
A_n(\bbU,\theta)\ \text{is real diagonalizable}
\quad\Longleftrightarrow\quad
A_H^N(T(\theta)\bbU)\ \text{is real diagonalizable}.
\end{equation}
Consequently, the system is (strongly) 2D hyperbolic at $\bbU$ if and only if $A_H^N(T(\theta)\bbU)$ is real diagonalizable for every $\theta\in[0,2\pi)$.
\end{proposition}

\begin{proof}
Equation \eqref{eq:rotation-invariance-HSWME} shows that $A_n(\bbU,\theta)=T(\theta)^{-1}A_H^N(T(\theta)\bbU)\,T(\theta)$, hence the two matrices are similar.
Similar matrices are real diagonalizable simultaneously, which yields \eqref{eq:equiv-diag}.
The final statement follows by Definition \ref{def:2D-hyperbolicity}.
\end{proof}

We next recall the characterization of real diagonalizability of $A_H^N$ proved in \cite{Bauerle2025}.

\begin{theorem}[Theorem 4.1 in \cite{Bauerle2025}]\label{thm:bauerle-41}
    Suppose $\alpha_1 = 0$.
    If $\beta_1 = 0$, then  $A^{N}_{H}(\bbU)$ in \eqref{mat:AH} is real diagonalizable; if $\beta_1 \neq 0$, then $A^{N}_{H}(\bbU)$ in \eqref{mat:AH} has only real eigenvalues but is not real diagonalizable.
\end{theorem}

\begin{theorem}[Theorem 4.4 in \cite{Bauerle2025}]\label{thm:bauerle-44}
    Suppose $\alpha_1 \neq 0$.
    Then $A^{N}_{H}(\bbU)$ in \eqref{mat:AH} is real diagonalizable.
\end{theorem}

Combining Theorem \ref{thm:bauerle-41} and Theorem \ref{thm:bauerle-44}, we obtain the characterization
\begin{equation}\label{eq:AH-diag-characterization}
A_H^N(\bbU)\ \text{is real diagonalizable}
\quad\Longleftrightarrow\quad
\bigl(\alpha_1\neq 0\bigr)\ \text{or}\ \bigl(\alpha_1=\beta_1=0\bigr).
\end{equation}

This statement should not be confused with the strong hyperbolicity of the two-dimensional HSWME.
The latter requires the condition in \eqref{eq:AH-diag-characterization} to hold for every rotated pair $(\alpha_1', \beta_1')$ generated by $T(\theta)\bbU$.

Now we are ready to present the proof to Theorem \ref{thm:HSWME}:
\begin{proof}[Proof to Theorem \ref{thm:HSWME}]
We first analyze the weakly hyperbolicity.
Since $A_H^N(\bbU)$ has only real eigenvalues for all $\bbU$, $A_n(\bbU,\theta)$ has only real eigenvalues for all $\bbU$ and $\theta$ by Proposition \ref{prop:reduction-AH}.
This is exactly weakly 2D hyperbolicity in Definition \ref{def:2D-hyperbolicity}.

We then analyze the strong hyperbolicity.
Under the rotation $T(\theta)$, the pair $(\alpha_1,\beta_1)$ transforms as
\begin{equation}\label{eq:rotated-alpha1-beta1}
    \alpha_1'=\cos\theta\,\alpha_1+\sin\theta\,\beta_1,\qquad \beta_1'=-\sin\theta\,\alpha_1+\cos\theta\,\beta_1.    
\end{equation}
\begin{enumerate}
    \item If $(\alpha_1,\beta_1)=(0,0)$ then $(\alpha_1',\beta_1')=(0,0)$ for all $\theta$, and \eqref{eq:AH-diag-characterization} implies $A_H^N(T(\theta)\bbU)$ is real diagonalizable for all $\theta$.
    Hence $A_n(\bbU,\theta)$ is real diagonalizable for all $\theta$, so the system is strongly 2D hyperbolic.

    \item If $(\alpha_1,\beta_1)\neq(0,0)$, one can always choose $\theta_0$ such that $\alpha_1'(\theta_0)=\cos\theta_0\,\alpha_1+\sin\theta_0\,\beta_1=0$.
    From \eqref{eq:rotated-alpha1-beta1}, $(\alpha_1')^2 + (\beta_1')^2 = (\alpha_1)^2 + (\beta_1)^2$ and thus $\beta_1'(\theta_0)\neq 0$.
    By \eqref{eq:AH-diag-characterization}, $A_H^N(T(\theta_0)\bbU)$ is not real diagonalizable.
    Therefore $A_n(\bbU,\theta_0)$ is not real diagonalizable, and the system is not strongly 2D hyperbolic at $\bbU$.
\end{enumerate}

\end{proof}

\FloatBarrier

\section*{Declaration of generative AI and AI-assisted technologies\\in the manuscript preparation process}
During the preparation of this work, the authors used OpenAI Codex for language refinement, LaTeX editing, consistency checking, and assistance with organizing the revision and submission materials.
The authors reviewed and edited the output as needed and take full responsibility for the content of the published article.

\bibliographystyle{abbrv}
\bibliography{ShallowWater}

\end{document}